%% file: main.tex
\documentclass[11pt]{article}
\usepackage{amsmath,amssymb,amsthm,bbm,color}
\usepackage[margin=1in]{geometry}
\usepackage{graphicx}
\usepackage{float}
\usepackage{placeins}
\usepackage{xcolor}

\usepackage{setup}
\usepackage{caption}

\usepackage{tikz}

\newcommand{\dcx}{\preceq_{\mathrm{dcx}}}
\newcommand{\stl}{\preceq_{\mathrm{st}}}
\newcommand{\stg}{\succeq_{\mathrm{st}}}
\newcommand{\cvx}{\preceq_{\mathrm{cvx}}}
\newcommand{\eqd}{\stackrel{\rm d}{=}}

\usetikzlibrary{arrows.meta,calc, patterns}

\title{Rolling Conformal Prediction in Sequential Model Training}

\author{%
	Chen Cheng\thanks{The first two authors contributed equally and are listed in alphabetical order.}\,\,\thanks{Department of Statistics, University of Illinois Urbana-Champaign} \and Ruiting Liang\footnotemark[1]\,\,\thanks{Department of Statistics, Harvard University} \and  Rina Foygel Barber\thanks{Department of Statistics, University of Chicago}
}

\begin{document}

\maketitle

\input{abstract}

\input{introduction}

\input{main-results}

\input{stability}

\input{numerical}

\input{discussion}

\section*{Acknowledgments}
C.C. and R.F.B. were supported by the Office of Naval Research via grant N00014-24-1-2544. R.F.B. was additionally supported by the National Science Foundation via grant DMS-2023109. R.L. was supported by the NSF-Simons AI-Institute for the Sky (SkAI) via grants NSF AST-2421845 and Simons Foundation MPS-AI-00010513.

\paragraph{Disclosure of AI use:} GPT-5.5 Plus, GPT-5.6 Plus, and GPT-5.6 Pro models were used for producing numerical experiments, assistance in writing, suggesting literature, and auditing existing proofs; the authors take full responsibility for the contents of the paper.
\bibliography{bib}
\bibliographystyle{abbrvnat}

\newpage

\appendix

\input{proof-main-results}

\input{proof-stability}

\input{appendix_additional_results}

\input{appendix_lemmas}

\end{document}

%% file: abstract.tex
\begin{abstract}
We introduce \textit{Rolling Conformal Prediction} (rolling-CP), a distribution-free predictive inference method for the setting of sequential model training. Specifically, given a data stream $(X_1,Y_1),(X_2,Y_2),\dots$, at each time $n$ the trained model may depend on the observed history $\{(X_i,Y_i)\}_{i<n}$. This setting arises naturally in modern sequential training, including one-pass training over massive datasets and continual fine-tuning or test-time adaptation of language models during deployment.

Rolling-CP first calibrates each incoming observation against the current predictor and then rolls it into future training. In this way, we avoid the need for data splitting. Remarkably, although the models at times $n=1,2,\dots$ may have entirely different properties and accuracy levels, for exchangeable data it is nonetheless possible to establish a guarantee of marginal coverage, with a familiar universal factor-two guarantee (a worst case guarantee of $1-2\alpha$ coverage, as compared to the target level $1-\alpha$), without any assumptions of stability or any restrictions on the model training process. For i.i.d.\ data streams, we further prove high-probability training-conditional validity uniformly over time; under stability conditions, coverage guarantees sharpen towards $1-\alpha$.
Numerical experiments on sequential regression, multiclass SGD, and one-pass neural-network training further demonstrate the practical effectiveness of rolling-CP.
\end{abstract}

%% file: introduction.tex
\section{Introduction}
The success of many modern machine learning systems depends on training sophisticated models through iterative procedures that continually update the model as new data are incorporated, such as stochastic gradient descent (SGD). These systems are often trained on massive datasets (e.g., LAION-5B for vision-language models and RefinedWeb for large language models \citep{schuhmann2022laion,penedo2023refinedweb}).

This evolving training process creates a fundamental challenge for predictive uncertainty quantification, since standard procedures are typically calibrated around a single fitted model. Conformal prediction \citep{vovk2005algorithmic,shafer2008tutorial} provides a generic and flexible framework for converting the outputs of a black-box model into prediction sets, in such a way as to guarantee (or closely approximate) the marginal coverage property
\begin{align}\label{eqn:intro-vague-coverage}
\P(Y_{n+1} \in \mc{C}_n(X_{n+1})) \approx 1 - \alpha,
\end{align}
where $\mc{C}_n(X_{n+1})$ is a prediction interval or prediction set for the unseen response $Y_{n+1}$, trained on the past data $(X_1,Y_1),\dots,(X_n,Y_n)$.
However, standard implementations of this framework do not directly address the challenge posed by an evolving training procedure: as we will describe in more detail shortly, existing conformal prediction methods are not designed to work within the setting of a dynamic training process.

This leads to a natural question: \emph{Can we design a distribution-free predictive inference procedure that exploits the sequential training trajectory to reduce computational cost while retaining a validity guarantee similar to \eqref{eqn:intro-vague-coverage}?}
In this paper, we answer this question affirmatively by introducing \emph{Rolling Conformal Prediction}, abbreviated as rolling-CP.

\subsection{Setting and background}\label{sec:setting_and_background}
We consider a streaming sequence of data points $(X_1,Y_1),(X_2,Y_2),\dots\in\mc{X}\times\mc{Y}$, where $X_i$ denotes a feature while $Y_i$ is the response. At time $n$,  the goal is to construct a prediction set $\mc{C}_n(X_{n+1})\subseteq\mc{Y}$ that is likely to contain the unseen response $Y_{n+1}$. 

To streamline the presentation, from this point on we will generally write $Z_i=(X_i,Y_i)$ to denote the $i$-th data point (and $z=(x,y)$ for a generic data point), and will reinterpret $\mc{C}_n$ as a subset of $\mc{Z} = \mc{X}\times\mc{Y}$---that is, 
\[z\in\mc{C}_n \ \Longleftrightarrow \ y\in\mc{C}_n(x)\textnormal{ for $z=(x,y)$.} \]
In other words, our goal is to provide a prediction set $\mc{C}_n\subseteq\mc{Z}$ that is likely to contain the next data point $Z_{n+1}$, by training on the previously observed data $Z_1,\dots,Z_n$.

Throughout the paper, we will primarily consider two regimes: the \emph{exchangeable setting}, where the goal is to provide a coverage guarantee at time $n+1$ under the assumption that $Z_1,\dots,Z_{n+1}$ are exchangeable, and the \emph{i.i.d.\ setting}, where we make the strictly stronger assumption that these $n+1$ data points are i.i.d.\ (from some unknown distribution). Throughout, we assume that $\mc{Z}$ is a standard Borel space, to ensure regularity conditions (e.g., for working with conditional distributions).

\paragraph{Background: split conformal prediction.} We begin with a brief overview of the split conformal prediction method \citep{vovk2005algorithmic,shafer2008tutorial}. Suppose $Z_1,\dots,Z_{n+1}$ are exchangeable (with $Z_i = (X_i,Y_i)$ as above). 
To construct $\mc{C}_n$ with the \emph{split conformal prediction} method, we first split the training data into two subsets---a subset $I\subseteq[n]$ used for model training, and the remaining data points $[n]\setminus I$ used for calibration. Let $n_0=|I|$ and $n_1 = n-n_0 = |[n]\setminus I|$.
Let $s: \mc{Z} \times \mc{Z}^{n_0} \to \R$ be a \emph{conformal score function}, where $s(z;(z_1,\dots,z_{n_0}))\in\R$ is a score indicating the extent to which $z$ does or does not ``conform'' to the trends observed in the data $(z_1,\dots,z_{n_0})$, with larger values indicating that $z$ may appear to be an outlier. In the setting of supervised learning with a real-valued response (i.e., the data consists of $(X,Y)$ pairs, with $Y\in\mc{Y}=\R$), a canonical example is the residual score: given a regression algorithm $\mc{A}$ that inputs training data and returns a fitted model $\widehat{f}:\mc{X}\to\R$, we can choose the score function as \begin{equation}\label{eqn:residual_score}s((x,y);(z_1,\dots,z_{n_0})) = |y - \widehat{f}(x)|\textnormal{ where $\widehat{f} = \mc{A}(z_1,\dots,z_{n_0})$}.\end{equation} 
The prediction set is then given by
\[\mc{C}_n = \left\{z\in\mc{Z} : s(z ; (Z_j)_{j\in I}) \leq \widehat{q}\right\}, \]
where the conformal quantile $\widehat{q}$ is defined as
\[\widehat{q} = \textnormal{Quantile}_{(1-\alpha)(1+1/n_1)}\left( (S_i)_{i\in [n]\setminus I}\right), \textnormal{ where } S_i = s(Z_i ; (Z_j)_{j\in I}).\]
By definition of the quantile, we can equivalently define this prediction set as
\begin{equation}\label{eqn:splitCP}\textbf{Split-CP:} \qquad \mc{C}_n = \left\{z\in\mc{Z} \ :  \ \sum_{i\in [n]\setminus I}   \ind\brc{s(z ; (Z_j)_{j\in I})>S_i} <(1-\alpha)(n_1+1)\right\}, \end{equation}
(This alternative definition will be useful for comparison to our proposed method, below.)
Under exchangeability of the data points $Z_1,\dots,Z_{n+1}$, this method offers a marginal coverage guarantee, $\P\prn{Z_{n+1}\in\mc{C}_n} \geq 1-\alpha$. However, it requires splitting the data, which may lead to a wider prediction interval (since the fitted model or score may be less accurate due to reduced sample size).

\subsection{The rolling conformal prediction method}
We are now ready to introduce the rolling-CP method.
For any $i\geq 1$, let $s_i:\mc{Z} \times \mc{Z}^{i-1} \to \R$ denote an arbitrary score function. Let $Z_{< i} = (Z_1,\dots,Z_{i-1})$ denote the data observed before time $i$, so that at this time our trained model is encoded in the score function $s_i(\cdot;Z_{<i})$, with the convention that a large value of the score $s_i(z;Z_{< i})$ indicates that $z$ does not appear to follow the trends observed in the training data $Z_{< i}$.\footnote{We use the convention that $Z_{<1} = \varnothing$, since no data has been observed before time $1$. The score function $s_1(z;Z_{<1})$ therefore depends only on $z$, and does not have access to any training data.}
For instance, if we are training a regression model using stochastic gradient descent (SGD) or some other iterative method, we might choose the residual score~\eqref{eqn:residual_score}, where $\widehat{f}_{<i}$ denotes the trained model after the $(i-1)$-st iteration (i.e., $\widehat{f}_{<i}$ is the model that was trained sequentially on the first $i-1$ data points).

We now define the rolling conformal prediction set, as follows:
\begin{align}
    \textbf{Rolling-CP:} \qquad \mc{C}_n =  \brc{z\in\mc{Z} \ : \  \sum_{i=1}^n \ind \brc{s_i(z; Z_{<i}) >s_i(Z_i; Z_{<i})} <(1-\alpha)(n + 1)}. \label{eqn:iterative-conformal-set}
\end{align}
The key feature of rolling-CP is that each incoming observation is used first as a fresh calibration point for the previously-trained score function and then rolled into the data set used to construct future score functions. Thus, the $i$-th comparison evaluates the candidate point $z$ and the newly arrived point $Z_i$ using the same score function $s_i(\cdot;Z_{<i})$. Aggregating these rolling comparisons yields the rolling-CP prediction set.

We can compare this definition to the split conformal method. At a high level, by examining the equivalent definition of split-CP given in~\eqref{eqn:splitCP}, we can see that the constructions are similar: each candidate value $z$ is compared to each calibration point $(Z_i)_{i\in [n]\setminus I}$, by comparing their respective scores; each such comparison (i.e., $z$ compared against each $Z_i$) uses the same shared trained score function $s(\cdot;(Z_j)_{j\in I})$. In contrast, rolling-CP compares each candidate value $z$ to \emph{every} data point $Z_i$ for $i\in[n]$ (rather than only comparing to a portion of the data), and each comparison uses a different score function (namely, for comparing $z$ and $Z_i$, we use the score function $s_i(\cdot;Z_{<i})$ that is the result of training up to time $i$). In other words, in rolling-CP every data point $i\in[n]$ is both used as a calibration point, and is then subsequently incorporated into training as the model evolves. As we will see in our experiments, this allows rolling-CP to make more efficient use of limited available data.

\paragraph{Organization and preview of main results.}
Before presenting our main guarantees for rolling-CP, we first illustrate how this framework captures several representative sequential training workflows in Section~\ref{sec:workflow-example}, and review related lines of work in Section~\ref{sec:related-work}.
Then, in Section~\ref{sec:main-results}, we present our main result: under exchangeability, rolling-CP offers the coverage guarantee
\begin{align*}
    \P\prn{Z_{n+1}\in\mc{C}_n} \geq 1-2\alpha.
\end{align*}
In Section~\ref{sec:stability}, we show that  coverage can be further sharpened towards the nominal coverage level $1-\alpha$, under additional  stability assumptions. Numerical experiments in Section~\ref{sec:numerical-experiment} further complement our theoretical results. Finally, we conclude with a brief summary in Section~\ref{sec:discussion}. Some proofs and additional results are deferred to the Appendix.

\subsection{Examples of sequential training settings}\label{sec:workflow-example}

For our first example, we begin with the setting of stochastic gradient descent (SGD) for streaming data.

\begin{example}[SGD for sequential training of a multiclass classifier] \label{example:SGD-classifier}
Let $Z_i=(X_i,Y_i)\in\R^D\times[K]$, $i=1,\dots,n+1$, be exchangeable, where $Z_1,\dots,Z_n$ arrive sequentially as training examples and $Z_{n+1}$ is the next test point. Consider a parametric class-index map
\begin{align*}
    f_\Theta:\R^D\to\R^K, \qquad
    f_\Theta(x)=\prn{f_{\Theta,1}(x),\dots,f_{\Theta,K}(x)},
\end{align*}
where $\Theta$ collects the trainable parameters. We refer to $f_\Theta(x)$ as the logit vector. In the linear multi-index model, $f_\Theta(x)=\Theta^\top x$ for $\Theta=(\theta_1,\dots,\theta_K)\in\R^{D\times K}$, and multinomial logistic regression models $\P(Y_i=k\mid X_i=x)=\mathsf{softmax}\prn{{\Theta^\star}^\top x}_k$, where $\mathsf{softmax}(v)_k:=\exp(v_k)/\sum_{\ell=1}^K\exp(v_\ell)$.

More generally, $f_\Theta$ may be a neural network, with $\Theta$ containing all network parameters. This formulation includes image-classification benchmarks such as CIFAR-10 ($K=10$) and CIFAR-100 ($K=100$) \citep{krizhevsky2009learning}, as well as the $1,000$-class ImageNet classification task \citep{russakovsky2015imagenet}. Starting from $\what{\Theta}_0$, suppose each arriving training example produces one SGD update,
\begin{align*}
    \what{\Theta}_i
    &=
    \what{\Theta}_{i-1}
    -
    \eta_i
    \left.
    \nabla_\Theta \ell(\Theta;Z_i)
    \right|_{\Theta=\what{\Theta}_{i-1}}, 
\end{align*}
where $\ell(\Theta;(x,y))=-\log\mathsf{softmax}\prn{f_\Theta(x)}_y$ is the cross-entropy loss. For $z=(x,y)$, one common choice are the cross-entropy score~\cite[Section~6.2]{goodfellow2016deep},
\begin{align*}
    s_i^{\mathsf{Ent}}(z;Z_{<i})
    =
    -\log\mathsf{softmax}\prn{f_{\what{\Theta}_{i-1}}(x)}_y.
\end{align*}
Another common choice is a running average of the logit-margin score,
\begin{align*}
    s_i^{\mathsf{Mar}}(z;Z_{<i})
    =
    \frac{1}{i\wedge T}
    \sum_{j=(i-T)_+}^{i-1}
    \brc{
        \max_{k\in[K]\setminus\{y\}} f_{\what{\Theta}_j,k}(x)
        -
        f_{\what{\Theta}_j,y}(x)
    }.
\end{align*}
The margin loss in multiclass formulation is standard~\cite{crammer2001algorithmic}, and we adopt the suffix averaging method of the most recent iterates as is common in stochastic optimization~\cite{rakhlin2012making}. Our experiments in
Sections~\ref{sec:numerical-sgd-logistic} and~\ref{sec:numerical-mnist} will study concrete instances of this setup for multi-index logistic regression and one-pass neural-network training.
\end{example}

Example~\ref{example:SGD-classifier} considers a fixed parametric model updated by SGD. The same ``rolling'' structure also appears naturally in modern LLM training, where the predictive pipeline may undergo more substantial changes.

While SGD is a very natural example, this type of sequential training can arise more broadly, in any training procedure that evolves over a data stream. This type of setting is especially relevant at modern scale, where models such as GPT-3 and Chinchilla operate in a single-pass or few-pass regime over their dominant data sources rather than over many repeated epochs \citep{brown2020language,hoffmann2022training}.

\vspace{0.3em}

\begin{example}[Continually updated language-model pipeline] \label{example:continual-llm} Let $Z_i=(X_i,Y_i)$ be exchangeable prompt--response pairs for $i \in [n+1]$. Before observing $Z_i$, the current language-model pipeline $M_{i-1}$ assigns a conditional distribution $p_{i-1}(\cdot\mid x)$ over a fixed response space $\mc{Y}$ (e.g., $\mc{Y}$ may be the set of candidate responses). A natural choice of the conformal score is the length-normalized negative log-likelihood~\cite{malinin2021uncertainty},
\begin{align*}
    s_i^{\mathsf{NLL}}((x,y);Z_{<i})
    =
    -\frac{1}{|y|}
    \sum_{t=1}^{|y|}
    \log p_{i-1}\prn{y_t\mid x,y_{<t}}.
\end{align*}
After $Z_i$ is scored, it may be used to update the pipeline. This update can be arbitrary (e.g., fine-tuning, replacing the model by a new checkpoint, or modifying its prompting components), as long as $M_{i-1}$ is determined from $Z_{<i}$. 
\end{example}

We note that both of these examples separate the statistical structure of the data stream from the evolution of the predictive model: the data stream itself can be exchangeable, while the latter may depend arbitrarily on the observed history.

\subsection{Related work}\label{sec:related-work}
In Section~\ref{sec:setting_and_background} above, we presented background on the split conformal method, which offers a marginal coverage guarantee in the setting of exchangeable data. Here we briefly mention several other methods within the conformal prediction framework that offer similar guarantees---in particular, methods that avoid the statistical cost of data splitting.
One alternative is the \emph{full conformal prediction} method \citep{vovk2005algorithmic}, which avoids data splitting but with a much steeper computational cost. Given a score function $s:\mc{Z}\times\mc{Z}^{n+1}\to\R$, the prediction set is given by
\[\mc{C}_n = \left\{z\in\mc{Z} : s\big(z ; (Z_1,\dots,Z_n,z)\big) \leq  \textnormal{Quantile}_{(1-\alpha)(1+1/n)}\left( S_1^z,\dots,S_n^z\right)\right\}, \]
for scores $S_i^z  = s\big(Z_i ; (Z_1,\dots,Z_n,z)\big)$. 
This method is computationally extremely expensive in general, requiring re-training the model or score function for each possible value $z\in\mc{Z}$ of the test point $Z_{n+1}$ (or, in the setting of supervised learning, for each possible value $y\in\mc{Y}$---that is, for all points of the form $z=(X_{n+1},y)$); in practice, accurate approximations to the conformal prediction set require a fine discretization of the response space \citep{lei2018distribution}
and remain computationally expensive. Other methods lie in between the statistical efficiency of full conformal, and the computational efficiency of split conformal: for example, cross-validation type versions of conformal prediction, including the cross-conformal method \citep{vovk2015cross,vovk2018cross} and the jackknife+ and CV+ methods \citep{barber2021predictive}. However, none of these aforementioned methods are designed to be computationally efficient in the setting of sequential model training.

Several lines of work have extended conformal prediction to settings of streaming data or sequential model training.
One direction is online conformal prediction, initiated by the adaptive conformal inference (ACI) method \citep{gibbs2021adaptive} and subsequently extended in \citep{gibbs2024conformal,bhatnagar2023improved, angelopoulos2024online}, among others. The core idea of this line of work is to calibrate the score threshold at each time using past coverage errors, thereby accommodating distribution shifts and dynamically updated models, potentially even under adversarial settings. Under minimal assumptions, these methods typically target long-run empirical coverage,
\begin{align*}
\frac{1}{T}\sum_{t=1}^T \ind\{Z_{t+1}\in\mc{C}_t\}
\to 1-\alpha,
\end{align*}
rather than marginal coverage at a fixed prediction time as \eqref{eqn:intro-vague-coverage}; stronger per-time conclusions require additional stochastic stability conditions \citep{angelopoulos2024online}. 
Although the running-average criterion is natural from a control-theoretic perspective \citep{angelopoulos2023pid} and can be enforced under essentially arbitrary data streams, it uses little of the probabilistic structure of the data. In contrast, we focus on exchangeable data streams and exploit this additional structure to show that rolling-CP prediction sets constructed from sequentially updated models can still enjoy valid per-time coverage guarantees. (See also Appendix~\ref{appendix_additional_results}, where we extend our guarantees to show robustness to mild distribution drift.)

Another line of work that extends conformal prediction to the sequentially arriving data setting, aims to constructing prediction sets that remain valid uniformly over time; see, for example, \citep{gauthier2025values, scharfstein2026time}. This direction is closely connected to the broader literature on anytime-valid inference \citep{ramdas2020admissible, howard2021time} and $e$-processes \citep{vovk2025conformal, ramdas2025hypothesis}. We note that, although the construction of rolling-CP is motivated by sequential settings, our main validity guarantee (see Theorem~\ref{thm:exchangeable-marginal}) is stated for any fixed sample size $n$. Nevertheless, rolling-CP can be deployed naturally in a sequential manner: $n$ need not be specified in advance, and the prediction set can be constructed at any time using the data and models available up to that point. More importantly, as we show in Section~\ref{sec:cond-unif}, this sequential construction also yields a time-uniform validity guarantee, which connects our method to the line of work of anytime-valid inference.

Lastly, we point out that rolling-CP can be viewed as aggregating a sequence of models and their pairwise comparison information in an efficient way, especially in terms of reducing the computational costs of both model-training and inference.
As a comparison, note that under the same sequential setup, existing methods for merging uncertainty sets \citep{solari2022multi, gasparin2024merging, gasparin2024conformal}, and more generally for merging $p$-values \citep{vovk2020combining, vovk2022admissible, wang2024testing, gasparin2025combining}, provide an alternative route to combining information across models while achieving the same marginal validity guarantee as ours.
Specifically, at time $n+1$, we can consider constructing split conformal prediction sets formed with various splits of the $n$ data points, and then merging these sets; we discuss a comparison between our method and this type of approach in Appendix~\ref{appendix_additional_results}.

%% file: main-results.tex
\section{Main results} \label{sec:main-results}
In this section, we present our main results establishing the coverage properties of rolling-CP in the exchangeable and i.i.d.\ data settings.

The rolling-CP prediction set is constructed in a fundamentally different way than traditional conformal prediction, and hence requires novel technical analyses. In traditional full or split conformal prediction, all scores entering the conformal rank are computed using a common trained model (e.g., the trained score function $s(\cdot;(Z_j)_{j\in I})$ in the definition~\eqref{eqn:splitCP} of split conformal prediction); the exchangeability of the data then guarantees that the test point's score is exchangeable with the scores of the other data point, leading to a direct argument establishing predictive coverage \citep{vovk2015cross,lei2018distribution,romano2019conformalized}. In contrast, rolling-CP repeatedly updates the score function as observations arrive and compares the same test point with scores constructed at different stages---the score at step $i$ may depend on previously revealed observations, and the usual arguments are no longer applicable.

Surprisingly, our findings reveal that under no additional structural assumptions of the scores $s_i(\cdot; Z_{<i})$, rolling-CP is able to achieve $1-2\alpha$ coverage guarantees for general exchangeable data stream. In particular, the score function $s_i(\cdot;Z_{<i})$ may be chosen in an arbitrary data-dependent way, may depend on the entire past history $Z_{<i}$, and need not satisfy any assumptions such as stability or consistency.

\subsection{Marginal coverage guarantee}\label{sec:marginal-coverage}

Our first result establishes marginal coverage under the sole assumption of exchangeability.

\begin{theorem}[Marginal coverage guarantee]\label{thm:exchangeable-marginal}
Let $Z_1,\dots,Z_{n+1}$ be exchangeable and let $\alpha\in(0,1)$. For any sequence of score functions $s_i:\mc{Z}\times\mc{Z}^{i-1}\to\R$,
the rolling-CP prediction set $\mc{C}_n$ defined in~\eqref{eqn:iterative-conformal-set} satisfies
\begin{align}
\P\prn{Z_{n+1}\in\mc{C}_n}\geq 1-2\alpha.\label{eqn}
\end{align}
\end{theorem}

For intuition, we begin by presenting a proof of this result in a simpler case, assuming that the data are i.i.d.\ and the scores have no ties almost surely (i.e., $s_i(Z_i;Z_{<i}) \ne s_i(Z_{n+1};Z_{<i})$ almost surely).

\paragraph{Proof of Theorem~\ref{thm:exchangeable-marginal}, special case (i.i.d.\ data, no ties).}
Assume that $Z_1,\dots,Z_{n+1}\simiid P$ for some distribution $P$, and assume that 
\[\textnormal{The conditional distribution of $s_i(Z;Z_{<i})\mid Z_{<i}$ is continuous, almost surely, when $Z\sim P$,}\]
to ensure that there are no ties among scores.
Define events $\mc{A}_1,\dots,\mc{A}_n$ and random variables $p_1,\dots,p_n$ as
\[\mc{A}_i = \brc{s_i(Z_{n+1};Z_{<i})>s_i(Z_i;Z_{<i})},\qquad p_i =  \P\prn{\mc{A}_i \mid Z_{<i+1}}.\]
Conditional on observing $Z_{<i}$, the random variables $s_i(Z_{n+1};Z_{<i})$ and $s_i(Z_i;Z_{<i})$ are i.i.d.\ and continuously distributed.
Therefore, by the standard CDF transform, it follows that $p_i \mid Z_{<i} \sim\mathsf{Unif}(0,1)$. Moreover, observe that $p_1,\dots,p_{i-1}$ are functions of $Z_{<i}$; consequently, by induction, we have shown that 
\[p_1,\dots,p_n \simiid \mathsf{Unif}(0,1).\]
Moreover, since the event $\mc{A}_i$ depends only on $Z_1,\dots,Z_i,Z_{n+1}$, and is independent of $Z_{i+1},\dots,Z_n$, we also have
$p_i = \P\prn{\mc{A}_i \mid Z_{<i+1}} = \P\prn{\mc{A}_i \mid Z_{<n+1}}$.

Next we need a lemma:
\begin{lemma}\label{lem:pigeonhole}
    Let $\mc{A}_1,\dots,\mc{A}_n$ be any events. Suppose that there exists some random variable $W$ such that 
    $\P\prn{\mc{A}_1 \mid W},\dots,\P\prn{\mc{A}_n \mid W}\simiid \mathsf{Unif}(0,1)$.
    Then for any $\Delta \geq 0$,
    \[\P\prn{\sum_{i=1}^n \ind_{\mc{A}_i} \geq  n-\Delta} \leq \frac{2\Delta +1}{n+1}.\]
\end{lemma}
\noindent Finally, by definition of the rolling-CP prediction set,
we have
\[Z_{n+1}\not\in\mc{C}_n \ \Longleftrightarrow \ \sum_{i=1}^n\ind_{\mc{A}_i} \geq (1-\alpha)(n+1).\]
Therefore, applying the lemma (with $\Delta=\alpha(n+1) - 1$, and $W=Z_{<n+1}$) completes the proof.
\endproof\bigskip

To complete the proof of this special case, we now prove the lemma.

\paragraph{Proof of Lemma~\ref{lem:pigeonhole}.}
Let $p_i = \P\prn{\mc{A}_i\mid W}$.
Write $p_{(1)}\leq \dots \leq p_{(n)}$ for the order statistics of $p_1,\dots,p_n$, and let $\mc{A}_{(j)}$ denote the event associated with $p_{(j)}$. Fix any positive integer $k$. Then, by the pigeonhole principle, 
\[\textnormal{If }\sum_{i=1}^n \ind_{\mc{A}_i} \geq n- \Delta\textnormal{ then }\sum_{j=1}^k \ind_{\mc{A}_{(j)}} \geq k- \Delta.\]
Therefore,
\[\P\prn{\sum_{i=1}^n \ind_{\mc{A}_i} \geq n- \Delta}
\leq \P\prn{\sum_{j=1}^k \ind_{\mc{A}_{(j)}} \geq k -  \Delta}
= \P\prn{\sum_{j=1}^k \ind_{\mc{A}_{(j)}} \geq k -  \lfloor\Delta\rfloor}
= \E\brk{\P\prn{\sum_{j=1}^k \ind_{\mc{A}_{(j)}} \geq k -  \lfloor\Delta\rfloor \,\bigg|\, W}},\]
where the second step holds since the sum $\sum_{j=1}^k \ind_{\mc{A}_{(j)}}$ is integer-valued. Applying Markov's inequality (and assuming $k > \lfloor \Delta\rfloor$), we then have
\[\P\prn{\sum_{i=1}^n \ind_{\mc{A}_i} \geq n- \Delta}
\leq \E\brk{\frac{\E\brk{\sum_{j=1}^k \ind_{\mc{A}_{(j)}}  \,\Big|\, W}}{k-\lfloor \Delta\rfloor}}=\E\brk{\frac{\sum_{j=1}^k p_{(j)}}{k-\lfloor \Delta\rfloor}} = \frac{\sum_{j=1}^k \frac{j}{n+1}}{k-\lfloor \Delta\rfloor}= \frac{k(k+1)}{2(n+1)(k-\lfloor \Delta\rfloor)},\]
where the next-to-last step holds since $p_1,\dots,p_n\simiid \mathsf{Unif}(0,1)$ and so $\E\brk{p_{(j)}} = \frac{j}{n+1}$. The above calculation holds for any  integer $k$ with $\lfloor \Delta\rfloor < k \leq n$; choosing $k = 2\lfloor \Delta\rfloor + 1$ completes the proof.
\endproof

The above proof strategy works only in the setting of i.i.d.\ data; the key step, where we show that $p_i\sim\mathsf{Unif}(0,1)$, may no longer hold under the weaker exchangeability assumption, since $Z_i$ and $Z_{n+1}$ may be dependent. As we will see below in Section~\ref{sec:cvx}, the general case can be proved using a different argument, via distributional ordering. (See also Appendix~\ref{proof:exchangeable-marginal} for an alternative proof that more closely follows the arguments used in the special case.)

In practice, we expect the coverage of rolling-CP to be approximately the nominal level, $1-\alpha$. However, the factor of two in Theorem~\ref{thm:exchangeable-marginal} cannot be improved without further assumptions such as convergence of the algorithm or stability of the score functions. Indeed, even for an i.i.d.\ uniform data stream, the asymptotic coverage can attain any value in $(1-2\alpha,1]$.

\begin{proposition}[Tightness of the coverage guarantee]\label{prop:tightness}
Let $\alpha\in(0,1/2)$ and $\nu\in(0,2\alpha]$. Let $P$ be any nonatomic distribution on $\mc{Z}$ (i.e., $\P(Z=z)=0$ for all $z\in\mc{Z}$, when $Z\sim P$). Then there exists a sequence of score functions $s_i:\mc{Z}\times\mc{Z}^{i-1}\to\R$ such that, for data $Z_1,Z_2,\dots\simiid P$, the rolling-CP method satisfies
\begin{align*}
\lim_{n\to\infty}\P\prn{Z_{n+1}\in\mc{C}_n}=1-2\alpha+\nu.
\end{align*}
\end{proposition}
\noindent However, as we will see later on in Section~\ref{sec:stability}, under an additional stability assumption we can ensure that the coverage of rolling-CP is approximately $1-\alpha$ (the nominal level).

\subsection{Proof of Theorem~\ref{thm:exchangeable-marginal} via convex ordering} \label{sec:cvx}
In this section, we establish a result characterizing an important property of rolling-CP, which relates to distributional ordering. This result will in turn lead to a proof of our main coverage guarantee, Theorem~\ref{thm:exchangeable-marginal}.

Distributional orderings compare random variables through their expectations over prescribed classes of test functions. These orderings provide a unified language for stochastic dominance and variability comparisons, e.g.~\cite{shaked2007stochastic, muller2002comparison}. Convex order is also the probabilistic analogue of majorization and is closely connected to doubly stochastic transformations and Lorenz order~\cite{marshall2011inequalities}. We begin with the relevant definitions.

\begin{definition}[Distributional orderings of random variables]Let $X,Y\in\R$ be random variables. We define the stochastic order, the convex order, and the decreasing convex order, respectively, as
\begin{align*}
X\stl Y&\quad\Longleftrightarrow\quad\E[f(X)]\leq \E[f(Y)] & & \textnormal{for every nondecreasing function $f:\R\to\R$}; \\
X\cvx Y&\quad\Longleftrightarrow\quad\E[f(X)]\leq \E[f(Y)] & & \textnormal{for every convex function $f:\R\to\R$}; \\
X\dcx Y&\quad\Longleftrightarrow\quad\E[f(X)]\leq \E[f(Y)]&&\textnormal{for every nonincreasing convex function $f:\R\to\R$}.
\end{align*}
The comparisons are restricted to measurable functions for which both expectations are well-defined. In particular, we can compare the definitions as follows:
\begin{equation*}\label{eqn:st_implies_dcx}
\textnormal{If $X\succeq_{\mathrm{st}}Y$ then $X\dcx Y$,}
\end{equation*}
and
\begin{equation*}\label{eqn:cvx_implies_dcx}
\textnormal{If $X\preceq_{\mathrm{cvx}}Y$ then $X\dcx Y$.}
\end{equation*}
\end{definition}

In hypothesis testing, a random variable $p\in[0,1]$ is a valid p-value (for testing a particular null hypothesis) if its distribution satisfies $\P(p\leq \alpha)\leq \alpha$ for all $\alpha\in(0,1)$ (under the null); in other words, $p\in[0,1]$ is a valid p-value if it satisfies $p\succeq_{\mathrm{st}}\mathsf{Unif}(0,1)$.
More generally,
\citet{wang2024testing} introduced the notion of a p*-value, which is any random variable $p\in[0,1]$ satisfying $p\dcx \mathsf{Unif}(0,1)$. This class contains ordinary p-values, and in fact is the convex hull of (possibly dependent) ordinary p-values. In particular, it is closed under arbitrary convex combinations without dependence assumptions, i.e., any average of p-values is a p*-value. A key fact is that any p*-value is itself, up to a factor of two, a valid p-value \citep{wang2024testing}:
\begin{equation}\label{eqn:p*-factor-of-2}
\textnormal{If $p$ is a p*-value (i.e., $p\in[0,1]$ and $p\dcx\mathsf{Unif}(0,1)$), then $\P\prn{p\leq \alpha}\leq 2\alpha$ for all $\alpha\in(0,1)$.}\end{equation}
To see why this is true, let $f(t) = (2-t/\alpha)_+$, which is a nonincreasing and convex function; we then observe that $\ind\brk{t\leq \alpha}\leq f(t)$ for all $t\in[0,1]$. Since $p\dcx U$ for $U\sim\mathsf{Unif}(0,1)$, we then have $\P\prn{p\leq \alpha}\leq \E\brk{f(p)} \leq \E\brk{f(U)} = 2\alpha$.

With this background in place, we are ready to proceed with the main theorem. Define a random variable
\[p_{\mathsf{rolling}} = \frac{1 + \sum_{i=1}^n \ind\{s_i(Z_i; Z_{<i}) \geq s_i(Z_{n+1}; Z_{<i})\}}{n+1}.\]
By construction, we have
\[Z_{n+1}\in \mathcal{C}_n\quad \Longleftrightarrow \quad p_{\mathsf{rolling}} > \alpha.\]
This quantity is analogous to the \emph{conformal p-value} \citep{vovk2005algorithmic,vovk2013transductive}, which offers an equivalent formulation of the split or full conformal prediction set. However, in the setting of rolling-CP, $p_{\mathsf{rolling}}$ is not necessarily a valid p-value---it might not necessarily satisfy $p_{\mathsf{rolling}}\stg\mathsf{Unif}(0,1)$.
The next theorem shows that when the score is instead updated arbitrarily over time, the quantity $p_{\mathsf{rolling}}$ instead satisfies a weaker distributional ordering property.

\begin{theorem}[Convex ordering dominance for exchangeable data stream] \label{thm:cvx-dcx} 
Let $Z_1,\dots,Z_{n+1}$ be exchangeable. For any sequence of score functions $s_i:\mc{Z}\times\mc{Z}^{i-1}\to\R$,
the random variable $p_{\mathsf{rolling}}$ satisfies
\[ p_{\mathsf{rolling}}\dcx \mathsf{Unif}\left(\left\{\frac{1}{n+1},\dots,\frac{n}{n+1},1\right\}\right).\]
\end{theorem}
\noindent This result leads directly to a proof of our main coverage guarantee---essentially, Theorem~\ref{thm:exchangeable-marginal} holds because, even though $p_{\mathsf{rolling}}$ is not necessarily a p-value, it is guaranteed to be a p*-value.

\paragraph{Proof of Theorem~\ref{thm:exchangeable-marginal}.}
The result of Theorem~\ref{thm:cvx-dcx}, along with the fact that $\mathsf{Unif}\left(\left\{\frac{1}{n+1},\dots,\frac{n}{n+1},1\right\}\right)\stg\mathsf{Unif}(0,1)$, implies that
\[p_{\mathsf{rolling}}\dcx\mathsf{Unif}(0,1).\]
In the terminology of \citet{wang2024testing}, this means that $p_{\mathsf{rolling}}$ is a p*-value, and consequently by~\eqref{eqn:p*-factor-of-2} we have $\P\prn{p_{\mathsf{rolling}}\leq \alpha}\leq 2\alpha$, for any $\alpha\in(0,1)$. Since $Z_{n+1}\in\mc{C}_n$ if and only if $p_{\mathsf{rolling}}>\alpha$, by construction, this means that
\[\P\prn{Z_{n+1}\in\mc{C}_n} = \P\prn{p_{\mathsf{rolling}}>\alpha}\geq 1-2\alpha,\]
as desired.
\endproof

\begin{remark}
    The coverage guarantee of Theorem~\ref{thm:exchangeable-marginal} can in fact be strengthened slightly. Since Theorem~\ref{thm:cvx-dcx} shows that $p_{\mathsf{rolling}}\dcx \mathsf{Unif}\left(\left\{\frac{1}{n+1},\dots,\frac{n}{n+1},1\right\}\right)$ (which is slightly stronger than proving that $p_{\mathsf{rolling}}$ is a p*-value, due to the difference between the discrete uniform distribution and $\mathsf{Unif}(0,1)$), it can be shown that this implies 
    $\P\prn{p_{\mathsf{rolling}}\leq \alpha}\leq \left(2\alpha - \frac{1}{n+1}\right)_+$,
    and consequently the marginal coverage guarantee can be improved to $\P\prn{Z_{n+1}\in\mc{C}_n}\geq 1 - (2\alpha - \frac{1}{n+1})_+$.
\end{remark}

\subsection{Conditional and uniform coverage in the i.i.d.\ case} \label{sec:cond-unif}
For many applications, it is also important to understand the behavior of the realized conformal set after observing the data stream, beyond the marginal coverage guarantee established in Theorem~\ref{thm:exchangeable-marginal}. We therefore next turn to the training-conditional coverage properties of rolling-CP. Throughout this section, we will work in the setting of i.i.d.\ data, with $Z_1,Z_2,\dots\simiid P$ for an arbitrary distribution $P$ on $\mc{Z}$.

The training-conditional coverage of a prediction set method is characterized by the random variable
\[\P\prn{Z_{n+1}\in\mc{C}_n\mid Z_{<n+1}},\]
which is the probability of coverage conditional on the training data $Z_{<n+1}=(Z_1,\dots,Z_n)$. 
For background, training-conditional coverage guarantees were established for split-CP by \citet{vovk2012conditional}. For full-CP, on the other hand, \citet{bian2023training} show that training-conditional coverage may fail arbitrarily badly without further assumptions, while \citet{liang2025algorithmic} establish that it holds under an algorithmic stability assumption. 

Our next result shows that rolling-CP guarantees  training-conditional coverage at a level $\geq 1-2\alpha - \mathcal{O}_P(n^{-\half})$, in the setting of i.i.d.\ data.
\begin{theorem}[Conditional coverage for i.i.d.\ data] \label{thm:training-conditional}
    Let $Z_1, \cdots, Z_n, Z_{n+1} \simiid P$ and let $\alpha \in (0,1)$. Define the rolling-CP prediction set $\mc{C}_n$ as in~\eqref{eqn:iterative-conformal-set}, and let
    \[\alpha_P(Z_{< n+1}) = \P(Z_{n+1} \notin \mc{C}_{n} \mid Z_1, \dots, Z_n)\]
    as the training-conditional miscoverage rate of $\mc{C}_n$. Then, for any $\delta \in (0,1)$,
    \begin{equation*}
        \P \prn{\alpha_P(Z_{< n+1}) \le 2\alpha+ \sqrt{\frac{2 \log(1/\delta)}{\alpha^2(n+1)}} } \ge 1-\delta.
    \end{equation*}
    In particular, this implies that for any $\delta\in(0,1)$, 
    \begin{align} \label{eqn:master-conditional-bound-uniform}
        \P \prn{ \alpha_P(Z_{< n+1}) \le  2\alpha+ \sqrt{\frac{2 \log ((n+1)^2/\delta)}{\alpha^2(n+1)}} \textnormal{ \ holds for all $n\geq 0$} } \ge 1-\delta.
    \end{align}
\end{theorem}
This uniform bound shows that, without additional assumptions, the training-conditional coverage probability at \emph{every time point} will uniformly be $\gtrapprox 1-2\alpha$. In particular, a direct implication of this uniform guarantee is that the marginal coverage of rolling-CP continues to hold approximately even when the test time is chosen adaptively from the observed data, for instance, at a suitable stopping time.

%% file: stability.tex
\section{Coverage guarantees under stability conditions} \label{sec:stability}

In this section, we show that, under additional conditions, the coverage guarantees established in Section~\ref{sec:main-results} can be further refined towards the nominal level $1-\alpha$.
The deviation from the nominal level is controlled by an additional error term that quantifies the extent to which the score functions $s_i(z; Z_{<i})$ \emph{stabilize} over the time $i=1,2,\dots$ of the data stream. In other words, given a test point $z$, does the score $s_i(z;Z_{<i})$ converge as $i\to\infty$? 

Throughout this section, we will work in the setting of i.i.d.\ data, with $Z_i\simiid P$ for an arbitrary distribution $P$.
Informally, our assumptions below will require that
\begin{multline}\label{eqn:informal_stab}\textnormal{After an initial period of $m$ time points, the score functions have stabilized,}\\\textnormal{ with $s_i(\cdot;Z_{<i})$ approximately equal to $s_j(\cdot;Z_{<j})$ for $i,j\geq m$.}\end{multline}

In continual training and online learning, this type of condition is a pathwise analogue of algorithmic stability, which measures sensitivity to perturbing the training sample~\citep{bousquet2002stability}. In conformal prediction, leave-one-out stability sharpens coverage for the jackknife and jackknife+~\citep{steinberger2018conditional,barber2021predictive} and yields training-conditional guarantees for full conformal and jackknife+~\citep{liang2025algorithmic,amann2023assumption}, typically with assumptions that bound score differences of the type
\begin{align*}
    |s_i(Z;Z_{<i}) - s_{i-1}(Z;Z_{<i-1})|.
\end{align*}
That is, the score is assumed to remain mostly unchanged if we remove one data point from the i.i.d.\ stream. In contrast, in the rolling-CP setting, assumptions of the form~\eqref{eqn:informal_stab} instead operate along a continually updated sequence: after time $m$, the score must remain stable despite many subsequent observations---the fitted model stabilizes to some limit.

\subsection{Results under score comparison stability}
We begin by formalizing the notion of stability that we will need for our results.

\begin{assumption}[Score comparison stability]\label{assmp:score-stability} Let $Z_1,Z_2,\dots\simiid P$, and let $Z, Z' \simiid P$ be drawn independently. Fix integers $n \geq m \geq 0$, and a parameter $\nu\in[0,1]$. Let $s_i : \mc{Z}\times\mc{Z}^{i-1}\to\R$ denote score functions, for $i\geq 1$, as before.
We assume that
    \begin{align} \label{eqn:score-stability-cond}
        \frac{1}{n-m+1} \sum_{i=m}^n \P\prn{\ind \brc{s_i(Z; Z_{<i}) > s_i(Z'; Z_{<i})} \neq \ind \brc{s_\star(Z; Z_{<m}) > s_\star(Z'; Z_{<m})}} \leq \nu,
    \end{align}
for some score function $s_\star : \mc{Z}\times\mc{Z}^{m-1}\to\R$.
\end{assumption}
\noindent Here, the score function $s_\star$ may simply be taken to be $s_m$, or it may be some limit point of the sequence. Essentially, this assumption quantifies the informal condition~\eqref{eqn:informal_stab} by requiring pairwise agreement (most of the time) between the score functions: if data points $Z,Z'$ are ranked in a particular order by the score function $s_i(\cdot;Z_{<i})$, this assumption ensures that they will likely also be ranked in the same order by $s_j(\cdot;Z_{<j})$ (where $i,j \geq m$).

\begin{theorem}[Marginal coverage under score comparison stability] \label{thm:stability-score} 
Let $Z_1, \cdots, Z_n, Z_{n+1} \simiid P$ and let $\alpha \in (0,1)$. Define the rolling-CP prediction set $\mc{C}_n$ as in~\eqref{eqn:iterative-conformal-set}, and suppose Assumption~\ref{assmp:score-stability} holds for some $m\in[n]$ and some $\nu\in[0,1]$. Then 
\begin{subequations}
\begin{align} \label{eqn:score-stability-coverage}
    \P(Z_{n+1} \in \mc{C}_{n})   \geq 1-\alpha\cdot \frac{n+1}{n-m+2} -2\sqrt{\nu}.
\end{align}
In addition, if the score function $s_\star$ does not lead to ties, almost surely (i.e., $\P(s_\star(Z;Z_{<m})=s_\star(Z';Z_{<m}))=0$, for data drawn i.i.d.\ from $P$), then we also have an analogous upper bound,
\begin{align}
    \P(Z_{n+1} \in \mc{C}_{n})   \leq 1-\alpha +2\sqrt{\nu} + \frac{(1-\alpha)(m-1)  + 1}{n-m+2}.
\end{align}
\end{subequations}
\end{theorem}
\noindent In other words, as long as the stability condition holds (for some $m\ll n$ and some $\nu\approx 0$), coverage is guaranteed to hold at approximately the nominal level $1-\alpha$.

We note that, under the score comparison stability assumption, the training-conditional coverage of rolling-CP can likewise be improved, to a rate $\geq 1-\alpha - \mc{O}_P(n^{-1/2})$ (as compared to $1-2\alpha - \mc{O}_P(n^{-1/2})$, as in Theorem~\ref{thm:training-conditional}). 
We defer these additional results to Appendix~\ref{appendix_additional_results}.

\subsection{Comparing to quantitative algorithmic stability conditions}
When studying stability of an algorithm, it is natural to directly compare the differences in fitted model values or score values, as a quantitative measure of stability---that is, to directly require that the values of the scores are approximately equal, $|s_i(Z;Z_{<i}) - s_j(Z;Z_{<j})|\approx 0$, rather than a condition based on pairwise comparisons, as in Assumption~\ref{assmp:score-stability}. This type of quantitative condition is more common in the literature (e.g.,~\citep{celisse2016stability}, as well as the examples mentioned at the the beginning of Section~\ref{sec:stability}). 
Here, we verify that our comparison-based stability condition is often weaker. 

We begin by stating a quantitative stability condition in the setting of sequential score training.
\begin{assumption}[$L_q$ score stability]\label{assmp:score-Lq-stability} In the setting and notation of Assumption~\ref{assmp:score-stability}, fix also any $q\in[1,\infty]$, and some parameter $\wt{\nu}_q\geq 0$. We assume that
    \begin{align} \label{eqn:score-Lq-stability-cond-with-sstar}
        \frac{1}{n-m+1} \sum_{i=m}^n  \E\brk{\norm{s_i(\cdot; Z_{<i}) - s_\star(\cdot; Z_{<m})}_{L_q}^q} \leq \wt{\nu}_q^q
    \end{align}
for $q < \infty$ or $ \max_{m \leq i \leq n} \norm{s_i(\cdot; Z_{<i}) - s_\star(\cdot; Z_{<m})}_{L_\infty} \leq \wt{\nu}_\infty$ almost surely with respect to $Z_{<n}$ for $q=\infty$.
\end{assumption}
\noindent (Here the norm $\norm{s_i(\cdot; Z_{<i}) - s_\star(\cdot; Z_{<m_\star})}_{L_q}$ is computed with respect to the distribution of $s_i(Z; Z_{<i}) - s_\star(Z; Z_{<m})$ induced by $Z\sim P$, and conditional on $Z_{<i},Z_{<m}$.)

We now verify that, under an additional bounded-density condition, $L_q$ score stability (Assumption~\ref{assmp:score-Lq-stability}) implies score comparison stability (Assumption~\ref{assmp:score-stability}).
\begin{proposition}[Connection between stability assumptions]\label{prop:convergence-Lq-implications}
Suppose Assumption~\ref{assmp:score-Lq-stability} holds for some $q\in[1,\infty]$ and $\wt{\nu}_q\geq 0$ and some choice of score function $s_\star:\mc{Z}\times\mc{Z}^{m-1}\to\R$.
If the conditional distribution of $s_\star(Z;Z_{<m})\mid Z_{<m}$ has density (with respect to Lebesgue measure) bounded by a constant $B$, almost surely, then
Assumption~\ref{assmp:score-stability} holds with parameter
\[\nu = \min\left\{ 1, 4(2B\widetilde{\nu}_q)^{q/(q+1)}\right\} .\]
\end{proposition}

If the bounded-density assumption does not hold, however, then this type of implication will no longer be meaningful (i.e., when $B$ can be arbitrarily large, this bound becomes vacuous due to $\nu=1$). More specifically, without a bounded density, it may be the case that an arbitrarily small perturbation of the score functions leads to a swap in the relative ranking of data points $Z,Z'$ with non-negligible probability. In this type of setting, the quantitative stability condition of Assumption~\ref{assmp:score-Lq-stability} can nonetheless be used to provide stronger coverage guarantees, if we slightly inflate the prediction set (similar to the results of \citep{barber2021predictive,liang2025algorithmic} for the jackknife+ and for full conformal); we omit the details for brevity.

%% file: numerical.tex
\section{Numerical experiments}\label{sec:numerical-experiment}

We complement our theoretical results with numerical experiments in three settings of increasing complexity: (i) Ordinary Least Squares (OLS) estimator for linear regression, (ii) Stochastic Gradient Descent (SGD) for logistic regression, and (iii) One-pass training on MNIST. The latter two experiments instantiate the sequential multiclass training setup of Example~\ref{example:SGD-classifier}, first for a linear multi-index model and then for neural-network logits. Together, the three examples range from analytically tractable estimators, to stochastic and large-scale iterative procedures, with score functions that may depend on the evolving fitted model.\footnote{Code for reproducing all experiments in this section is available at \url{https://github.com/Moriartycc/rolling-conformal}}

For each experiment, we construct rolling-CP sets at several nominal coverage levels $1-\alpha$, using data-dependent score functions. We repeat the complete experiment over independent draws of both the dataset and, when applicable, the algorithmic randomization (e.g., initialization). In each experiment, we are particularly interested in the empirical marginal coverage of rolling-CP (as compared to the nominal level $1-\alpha$), which we evaluate using an independent hold-out set.

\subsection{OLS for linear regression} 
We begin with the simple linear model where $X_i \simiid \normal(0, I_d)$, with the ground-truth $\theta^\star = e_1 \in \R^d$ and observations $Y_i=X_i^\top \theta^\star + \varepsilon_i$, where $\varepsilon_i \simiid \normal(0, \sigma^2)$ are independent of the features. The data stream is then formed by $Z_i=(X_i, Y_i)$. We consider the sequence of minimum-norm OLS estimators. Let $\what{\theta}_0 = 0$ and for $i \in [n]$,
\begin{align*}
    \what{\theta}_i = \argmin \limits_{\theta \in \R^d} \frac{1}{2i} \sum_{j=1}^i \norm{X_j^\top \theta - Y_j}_2^2, \qquad \textrm{s.t. }\, \norm{\theta}_2 \, \textrm{ is minimized}.
\end{align*}
We can also write its matrix form concisely as $\what{\theta}_i = X_{1:i}^\dagger Y_{1:i}$ where $X_{1:i} = (X_1, \dots, X_i)^\top \in \R^{i \times d}$ and $Y_{1:i}=(Y_1, \dots, Y_i)^\top \in \R^i$. In this example, for $z=(x,y)$, we fix the sequence of scores as
\begin{align*}
    s_i(z; Z_{<i}) = \frac{1}{2}(x^\top \what{\theta}_{i-1} - y)^2.
\end{align*}

We can then use this paradigm to study predictive uncertainty across the interpolation threshold in a continual-learning regime by letting the data stream size $i$ increase from the overparameterized regime $i<d$, through the interpolation threshold $i\approx d$ and into the underparameterized regime $i>d$. In this setting, prediction risk exhibits double descent and instability near interpolation \citep{belkin2019reconciling,hastie2022surprises}, while exact interpolation can still yield good generalization via benign overfitting \citep{bartlett2020benign, cheng2024dimension}, and may even be necessary for optimal risk in some high-dimensional models \citep{cheng2022memorize, cheng26memorization}. 

\begin{figure}[!htbp]
    \centering
    \begin{minipage}[t]{0.485\textwidth}
        \centering
        \includegraphics[width=0.94\linewidth]{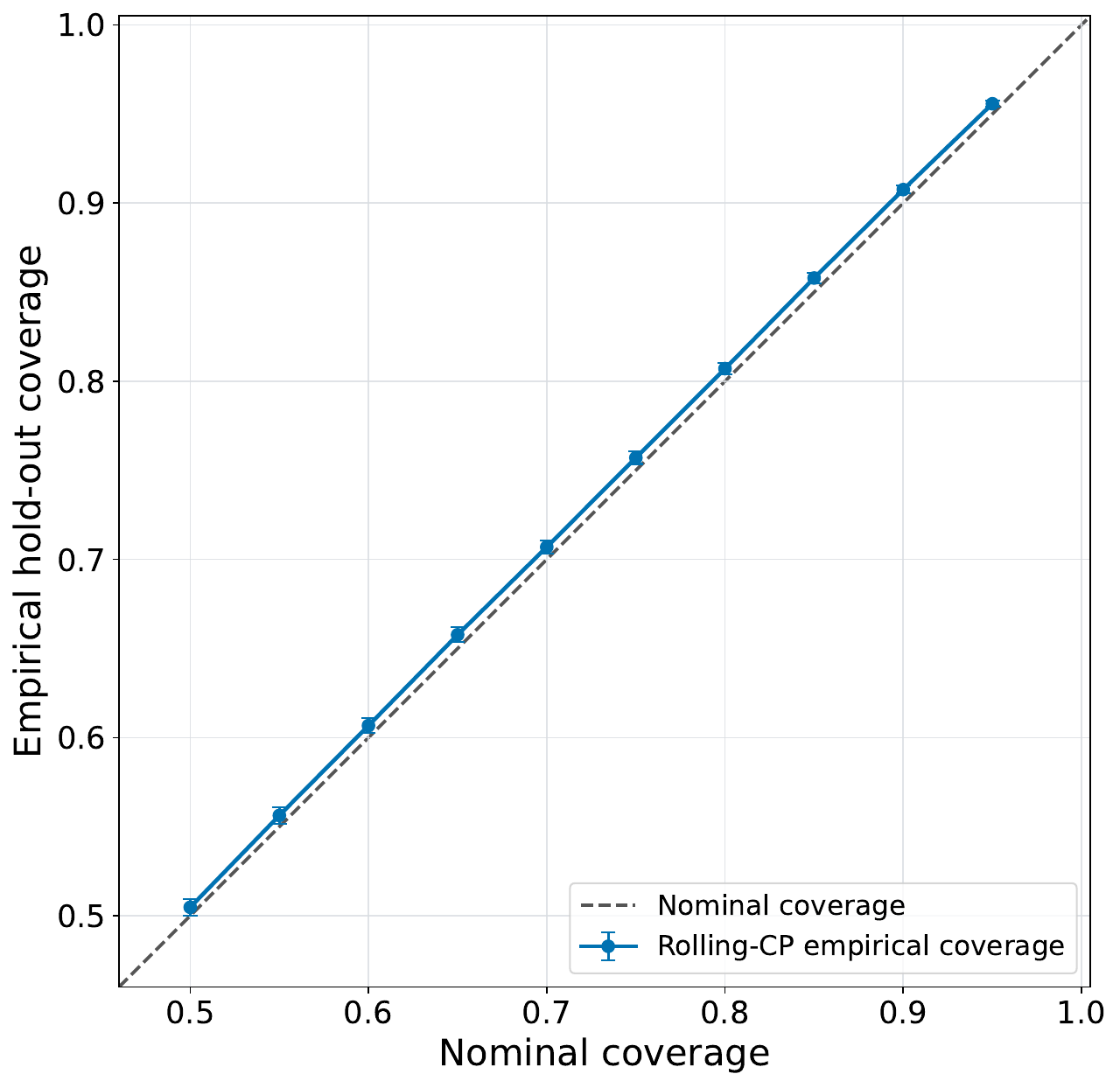}
        \vspace{0.3em}
        
        \parbox{\linewidth}{\footnotesize \textbf{(a)} For each nominal level, we report the mean empirical coverage at \(i=n\) across the \(M=100\) data streams.}
    \end{minipage}
    
    \vspace{0.8em}
    
    \begin{minipage}[t]{0.485\textwidth}
        \centering
        \includegraphics[width=1.05\linewidth]{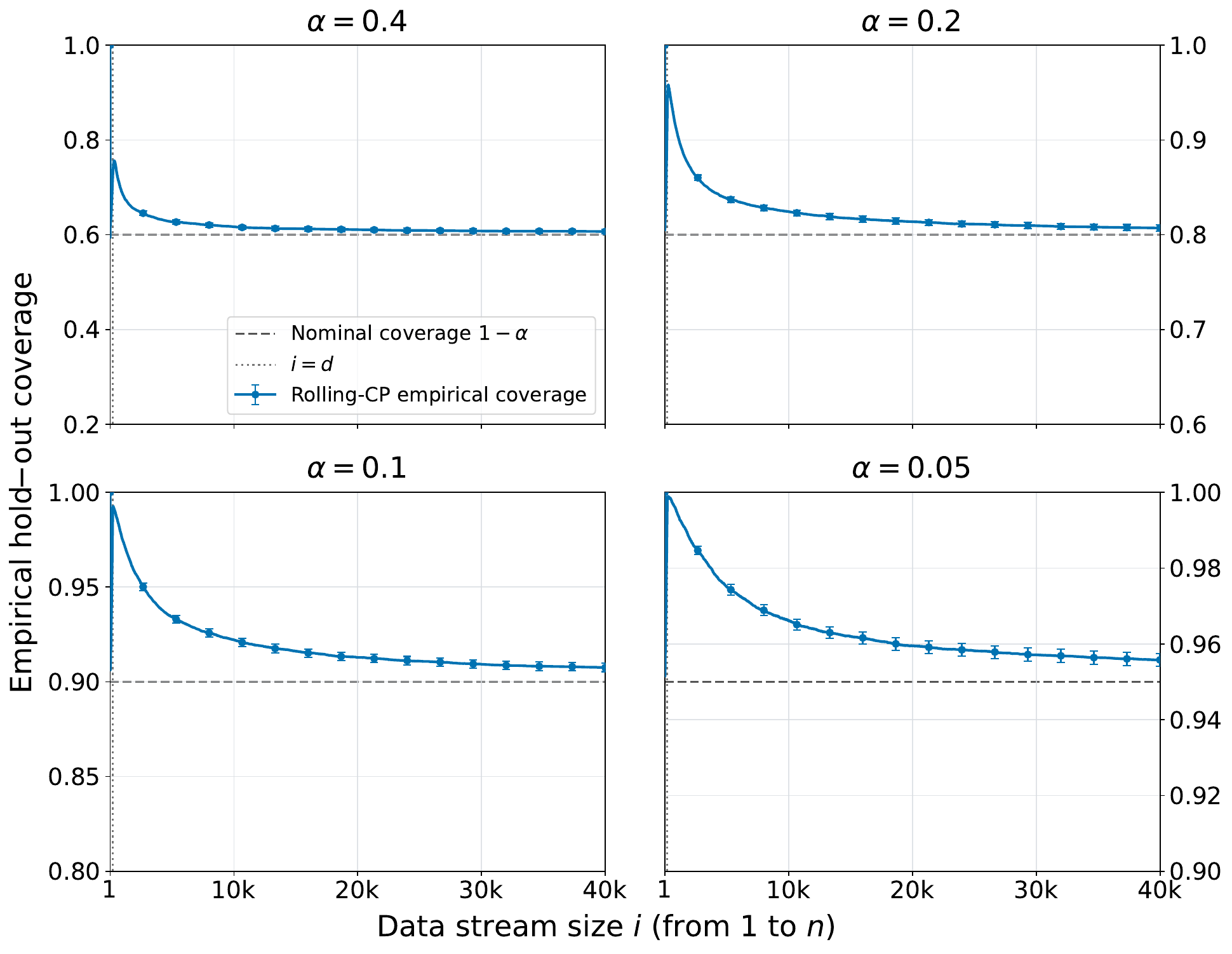}
        \vspace{0.3em}
        
        \parbox{\linewidth}{\footnotesize \textbf{(b)} For $\alpha\in\{0.4,0.2,0.1,0.05\}$, we report average empirical coverage along the data stream. Horizontal dashed lines are the nominal level $1-\alpha$, the vertical dotted line is the interpolation threshold $i=d$ (near the $y$-axis).}
    \end{minipage}%
    \hfill
    \begin{minipage}[t]{0.485\textwidth}
        \centering
        \includegraphics[width=1.05\linewidth]{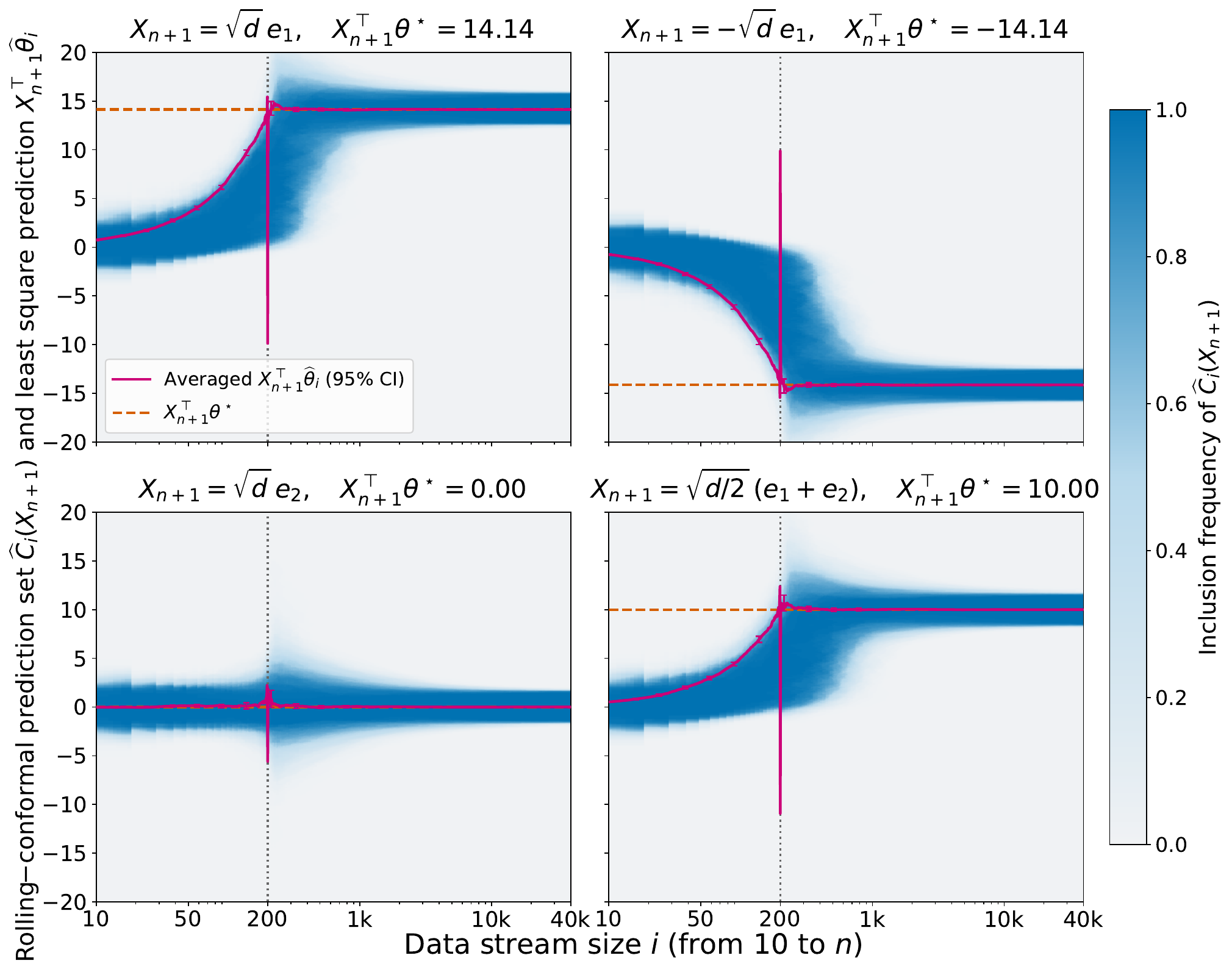}
        \vspace{0.3em}
        
        \parbox{\linewidth}{\footnotesize \textbf{(c)} The \textcolor{blue}{shaded area} indicates the inclusion frequency in the $90\%$ rolling-CP set, the \textcolor{magenta}{curve} is the averaged least-squares point prediction and the \textcolor{orange}{dashed line} is the oracle conditional mean.}
    \end{minipage}

    \vspace{1em}
    
    \caption{Numerical experiments of OLS for linear regression. \textbf{(a)} End-of-stream empirical marginal coverage. \textbf{(b)} Evolution of coverage across the overparameterized, interpolation, and underparameterized regimes. \textbf{(c)} Four fixed-feature prediction sets and least-squares point predictions. All panels use $d=200$, $n=40,000$, $\sigma=1$, and $M=100$ independent data streams with $95\%$ confidence intervals for the error bars.}
    \label{fig:linear-regression}
\end{figure}

\paragraph{Rolling conformal prediction applied to the sequential data stream.} 
Figure~\ref{fig:linear-regression} shows the performance of rolling-CP in this setting. We set $d=200$, $n=40,000$, and $\sigma=1$, and average the results over $M=100$ independent data streams. We compute the full min-norm least-squares trajectory and construct rolling-CP sets using the squared residual scores defined above. We estimate marginal coverage on an independent hold-out sample and use several fixed test-feature values to visualize how the prediction sets evolve. The three panels of Figure~\ref{fig:linear-regression} report the resulting coverage of rolling-CP, its evolution along the data stream, and the test-conditional prediction sets. Note that the empirical coverage is extremely close to the nominal level $1-\alpha$ (suggesting that this setting may exhibit stability properties).

    \paragraph{Rolling-CP vs.\ split-CP.} It is instructive to compare rolling-CP with split-CP, which is also natural in this setting. Recall the definition of split-CP~\eqref{eqn:splitCP}: at any time $n$, we use data points $I\subseteq[n]$ for model training, and the remaining points $[n]\setminus I$ as a calibration set. In contrast, rolling-CP does not require data splitting.

    Here we consider two settings. 
    In the first setting, we take $|I| = m$ to be a \emph{fixed} training set size. Then,
    \begin{itemize}
        \item For split-CP, we train on data points $I = \{1,\dots,m\}$; then, for $n\geq m$, at time $n+1$ we use data points $[n]\setminus I = \{m+1,\dots,n\}$ as a calibration set, and data point $n+1$ as a test point.
        \item For rolling-CP, we treat the first $m$ data points as a burn-in period; that is, at time $n+1$, we apply rolling-CP~\eqref{eqn:iterative-conformal-set} using only data points $i=m+1,\dots,n$ as our data stream (rather than data points $i=1,\dots,n$) for producing a prediction set for test point $n+1$.
    \end{itemize}
    In the second setting, we run split-CP with an \emph{increasing} training set size, $|I|=\lfloor n/2\rfloor$. Then,
    \begin{itemize}
        \item For split-CP, at time $n+1$, we train on data points $I = \{1,\dots,\floor n/2\rfloor\}$, and we use data points $[n]\setminus I = \{\lfloor n/2\floor+1,\dots,n\}$ as a calibration set, and data point $n+1$ as a test point.
        \item Rolling-CP is implemented exactly as in~\eqref{eqn:iterative-conformal-set}, with no burn-in period.
    \end{itemize}

\begin{figure}[!htbp]
    \centering
    \begin{minipage}[t]{0.485\textwidth}
        \centering
        \includegraphics[width=1.05\linewidth]{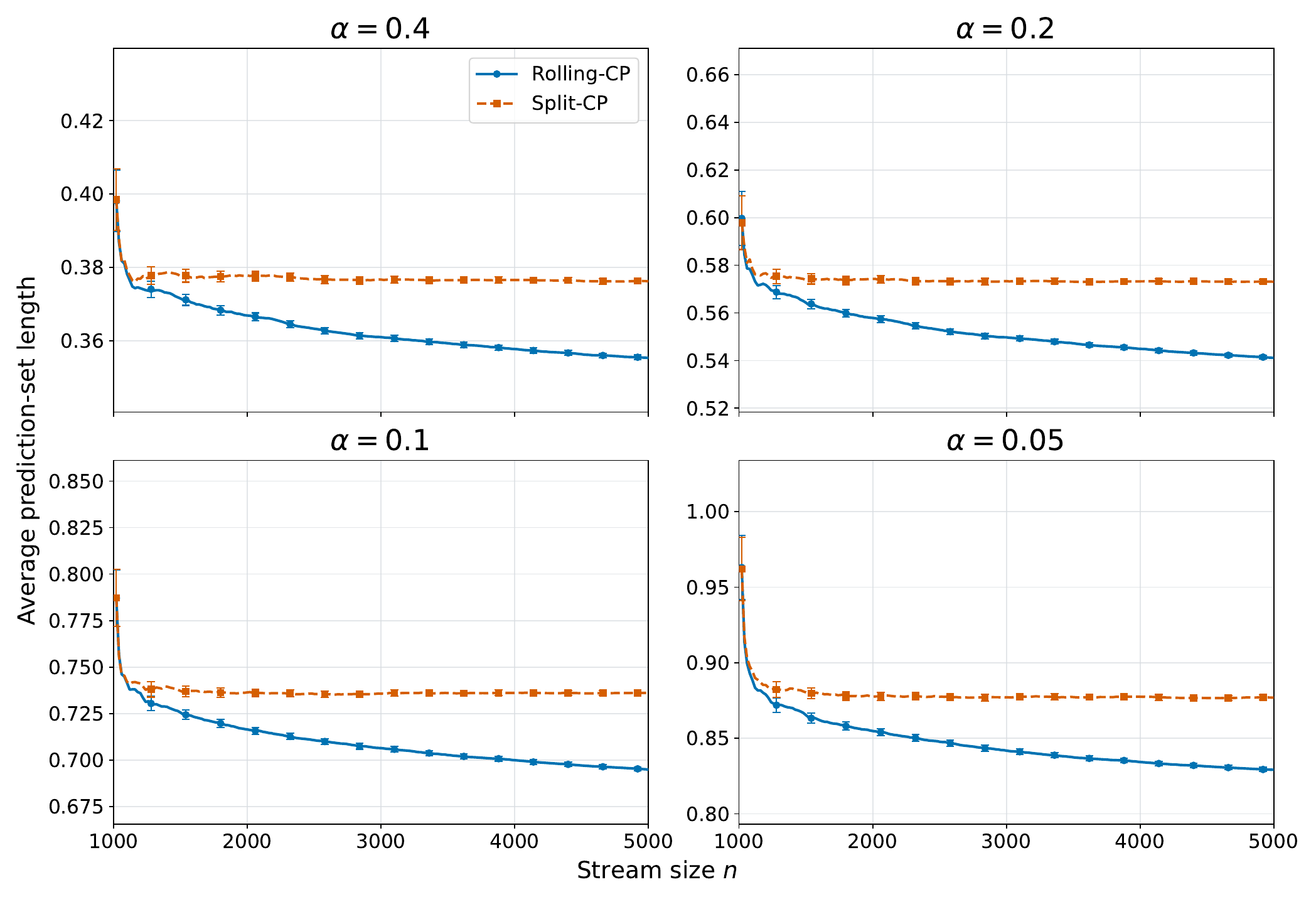}
        \vspace{-0.4em}

        {\footnotesize \textbf{(a)} Average length, fixed $m=1,000$.}
    \end{minipage}
    \hfill
    \begin{minipage}[t]{0.485\textwidth}
        \centering
        \includegraphics[width=1.05\linewidth]{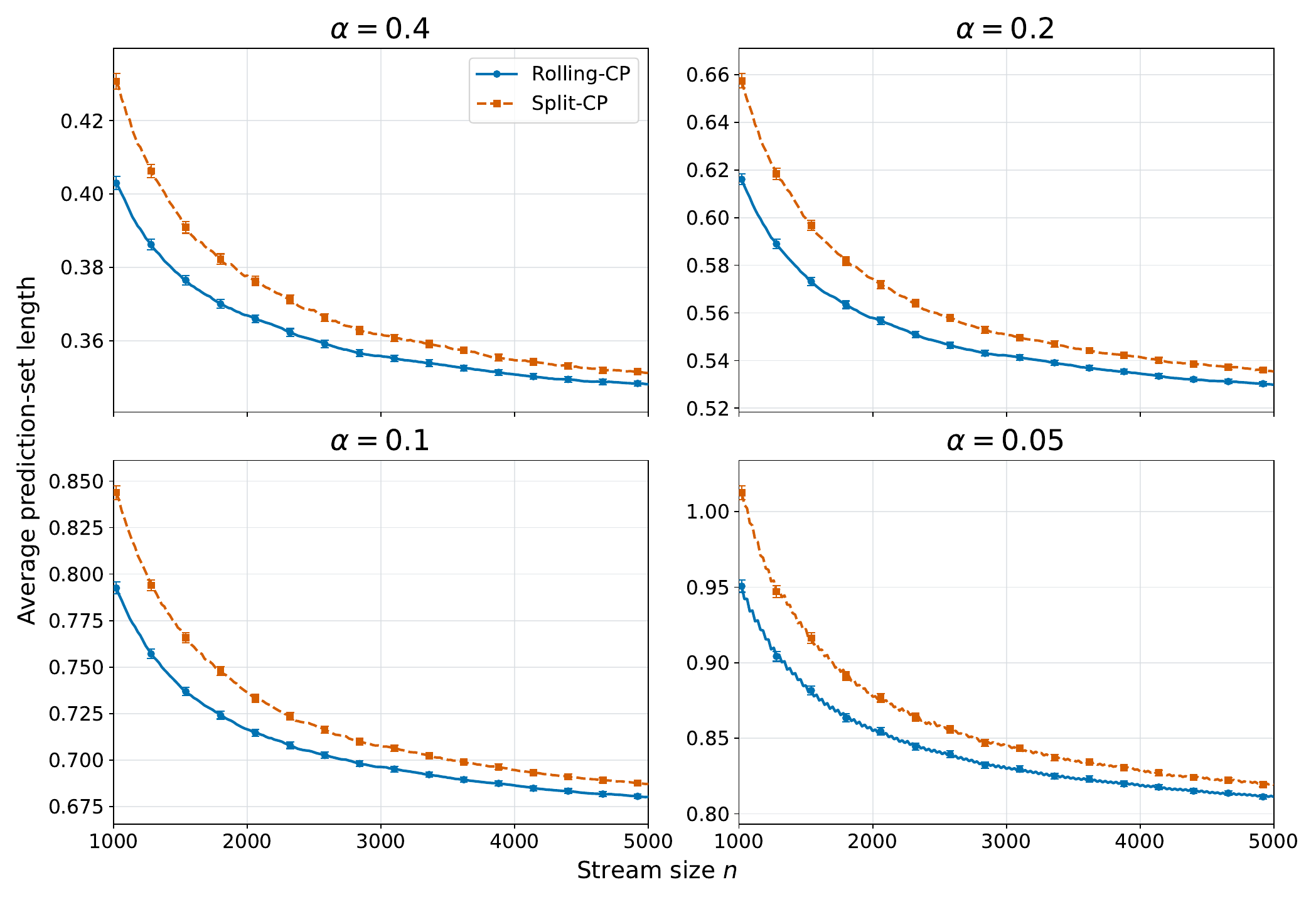}
        \vspace{-0.4em}

        {\footnotesize \textbf{(b)} Average length, $m=\lfloor n/2\rfloor$.}
    \end{minipage}

    \vspace{0.6em}
    \begin{minipage}[t]{0.485\textwidth}
        \centering
        \includegraphics[width=1.05\linewidth]{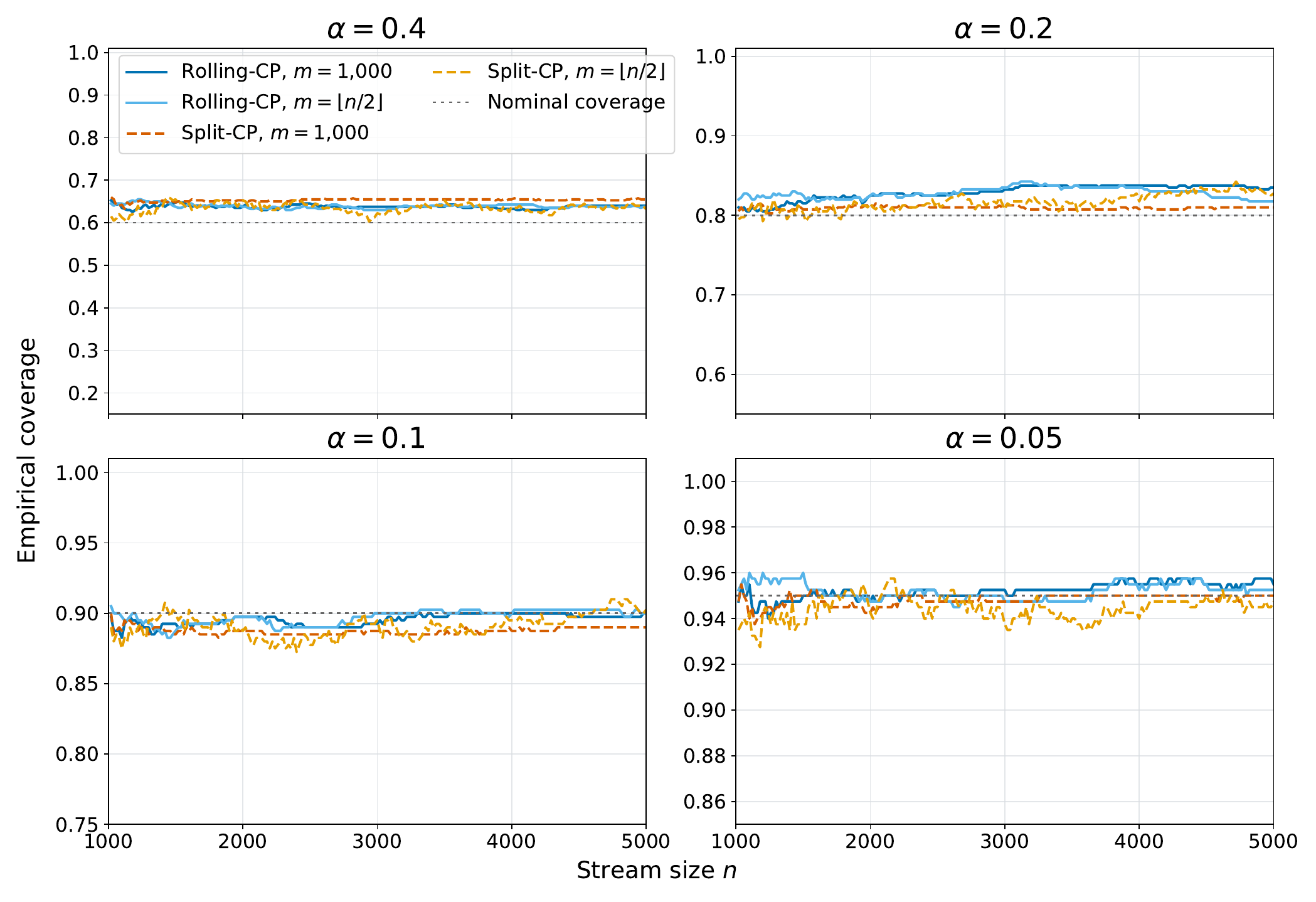}
        \vspace{-0.4em}

        {\footnotesize \textbf{(c)} Empirical coverage.}
    \end{minipage}
    \vspace{.6em}
    \caption{Comparison of rolling-CP and split-CP for minimum-norm OLS. \textbf{(a)} Evolution of the average prediction-set length with the fixed cutoff $m=1,000$. \textbf{(b)} Evolution of the average prediction-set length with the proportional cutoff $m=\lfloor n/2\rfloor$. \textbf{(c)} Evolution of empirical coverage for both methods and cutoff choices; horizontal dotted lines indicate the nominal coverage $1-\alpha$. All panels use $d=200$, terminal sample size $n=5,000$, $\sigma=0.2$, and $M=400$ independent data streams. Error bars in panels \textbf{(a)} and \textbf{(b)} are pointwise $95\%$ confidence intervals across trials.}
    \label{fig:linear-regression-vs-split}
\end{figure}

We present the results of this comparison in Figure~\ref{fig:linear-regression-vs-split}. We set $d=200$ and $\sigma=0.2$, take the terminal sample size to be $n=5000$, and evaluate both methods at $n=1020,1040,\ldots,5000$. In the first setup, we fix $m=1000$. In the second setup, we take $\lfloor n/2\rfloor$ training points at each $n$. For each configuration, we estimate marginal coverage and average prediction-set length over $M=400$ independent data streams, using the same data stream and independent test point for all methods within each trial. 

The three panels of Figure~\ref{fig:linear-regression-vs-split} report prediction-set length under the two setups and empirical coverage along the data stream. We observe that the prediction sets for rolling-CP are substantially narrower as $n$ increases (while coverage is similar for the two methods), highlighting the benefit of avoiding data splitting in order to use the data more efficiently.

\subsection{SGD for logistic regression}\label{sec:numerical-sgd-logistic} 
We next consider a $K=5$ class multi-index logistic-regression model. Consider the same feature distribution as in the linear regression example, $X_i \simiid \normal(0, I_d)$ and, conditional on $X_i=x$, generate $Y_i\in[K]$ according to
\begin{align*}
    \mathbb{P}(Y_i=k\mid X_i=x)=\mathsf{softmax}\prn{{\Theta^\star}^\top x}_k,
\end{align*}
where we take the following ground-truth indices $\Theta^\star = (\theta_1^\star, \dots, \theta_K^\star) \in \R^{d \times K}$:
\begin{align*}
\theta_1^\star=e_1,\qquad\theta_2^\star=e_2,\qquad\theta_3^\star=e_3,\qquad\theta_4^\star=\frac{1}{2}(e_1+e_2),\qquad\theta_5^\star=\frac{1}{2}(e_2+e_3+e_4+e_5)
\end{align*}
Starting from $\what{\Theta}_0=0$, where the columns of $\what{\Theta}_i \in\mathbb{R}^{d \times K}$ are the estimated class coefficients, we process each observation once using online SGD on $\R^{d \times K}$ with a polynomial decaying stepsize:
\begin{align*}
    \what{\Theta}_i=\what{\Theta}_{i-1}-\eta_i X_i\prn{\mathsf{softmax}\prn{\what{\Theta}_{i-1}^\top X_i}-e_{Y_i}}^\top, \qquad \eta_i=\frac{\eta_0}{(t_0+i)^\gamma},
\end{align*}
where $e_{Y_i}$ is the $Y_i$-th basis vector in $\R^K$, $\frac{1}{2} < \gamma \leq 1$ is the decaying exponent, chosen to satisfy the classical Robbins-Monro conditions~\cite{robbins1951stochastic}, and $\eta_0, t_0 > 0$ are initializing constants. For $z=(x,y) \in \R^d \times [K]$, we consider two sequences of scores: (i) the cross-entropy score $s_i^{\mathsf{Ent}}$ and (ii) the running average of logit-margin score $s_i^{\mathsf{Mar}}$ with a moving window size $T$:
\begin{align*}
    s_i^{\mathsf{Ent}}(z; Z_{<i}) = -\log \mathsf{softmax}\prn{\what{\Theta}_{i-1}^\top x}_y, \qquad s_i^{\mathsf{Mar}}(z; Z_{<i}) = \frac{1}{i \wedge T} \sum_{j=(i-T)_+}^{i-1} \brc{\max_{k\in[K]\setminus\{y\}} e_k^\top\what{\Theta}_j^\top x - e_y^\top\what{\Theta}_j^\top x}.
\end{align*}

This experiment tests rolling-CP in a genuinely iterative training regime. During a single pass through the data, the SGD iterations can fluctuate substantially before stabilizing to their limiting behavior. The validity of rolling-CP does not require the SGD iterates to have converged, nor does it require stability or consistency of the fitted classifier.

We report our numerical results in Figure~\ref{fig:SGD-logistic}, taking $d=10$, $n=10,000$, $\eta_0=1$, $t_0=10$, $T=100$ and averaging the results over $M=100$ independent data streams. For each $\gamma\in\{0.6,0.8,1\}$, we construct rolling-CP using both the cross-entropy and moving-average logit-margin scores defined above. We estimate marginal coverage on an independent hold-out sample and use four fixed test-feature values to visualize how the prediction sets evolve. The six panels of Figure~\ref{fig:SGD-logistic} report end-of-stream coverage, coverage along the data stream, and class-wise inclusion frequencies. Overall, we see that rolling-CP achieves approximately the nominal coverage level empirically, despite the fluctuations of the trained model in this setting.

\begin{figure}[!htbp]
\centering
\begin{minipage}[t]{0.45\textwidth}
\centering
\includegraphics[width=1.02\linewidth]{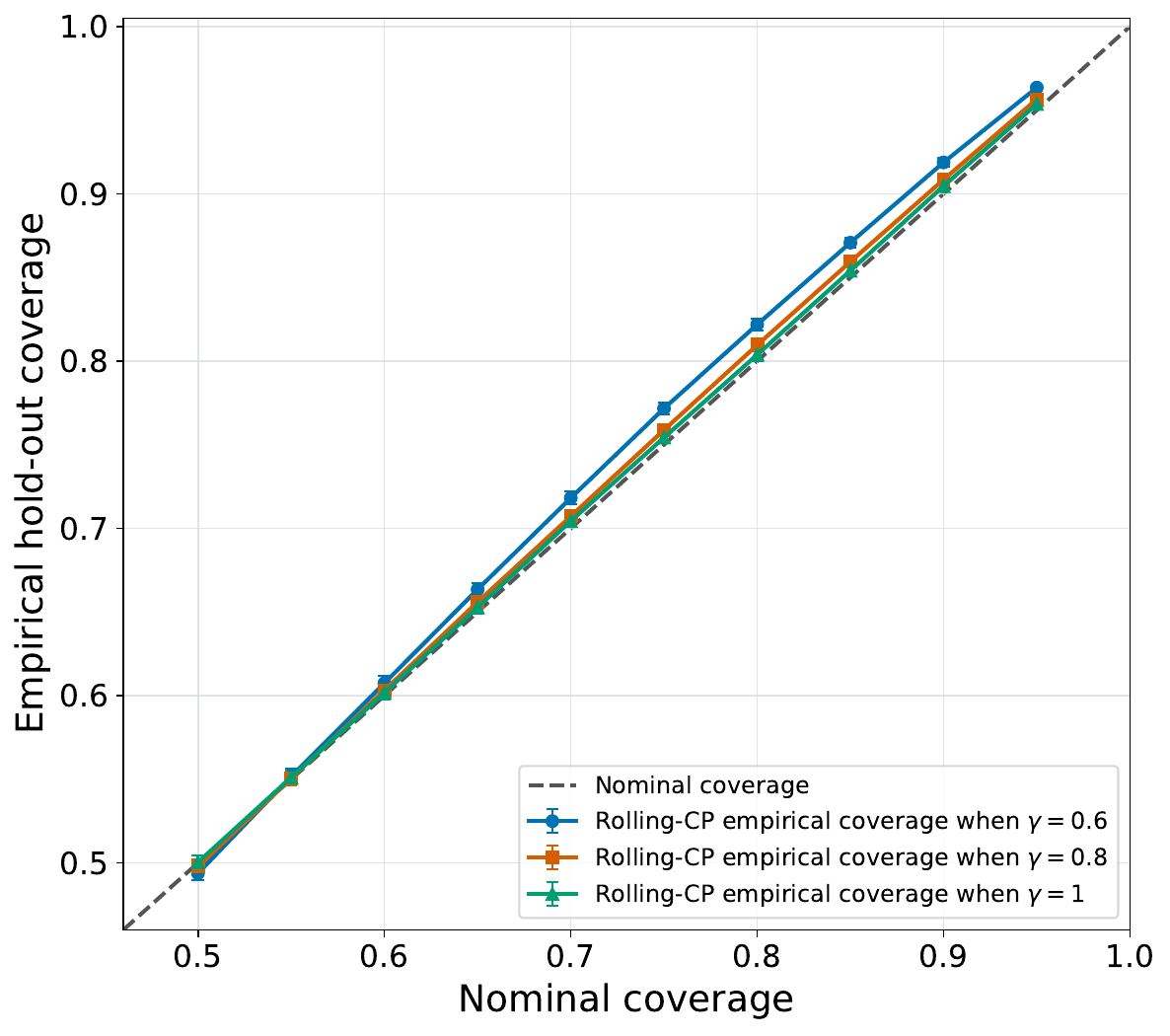}
\par\vspace{0.1em}
\parbox{\linewidth}{\footnotesize \textbf{(a)} End-of-stream held-out coverage versus the nominal level for $s_i^{\mathsf{Ent}}$ under choices of $\gamma$.}
\end{minipage}%
\hspace{0.03\textwidth}
\begin{minipage}[t]{0.45\textwidth}
\centering
\includegraphics[width=1.02\linewidth]{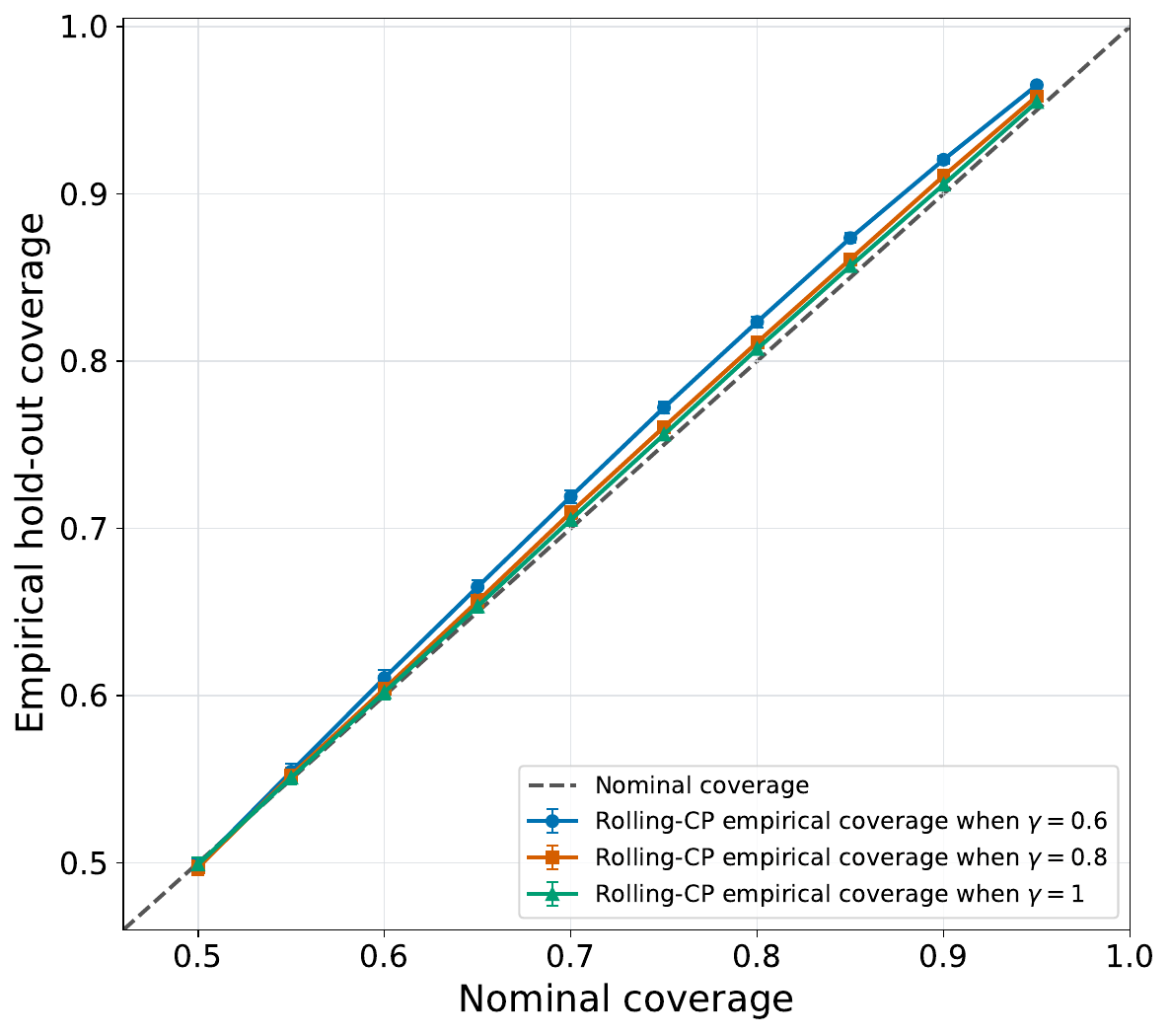}
\par\vspace{0.1em}
\parbox{\linewidth}{\footnotesize \textbf{(b)} End-of-stream held-out coverage versus the nominal level for $s_i^{\mathsf{Mar}}$ under choices of $\gamma$.}
\end{minipage}
\par\vspace{0.3em}

\begin{minipage}[t]{0.45\textwidth}
\centering
\includegraphics[width=1.02\linewidth]{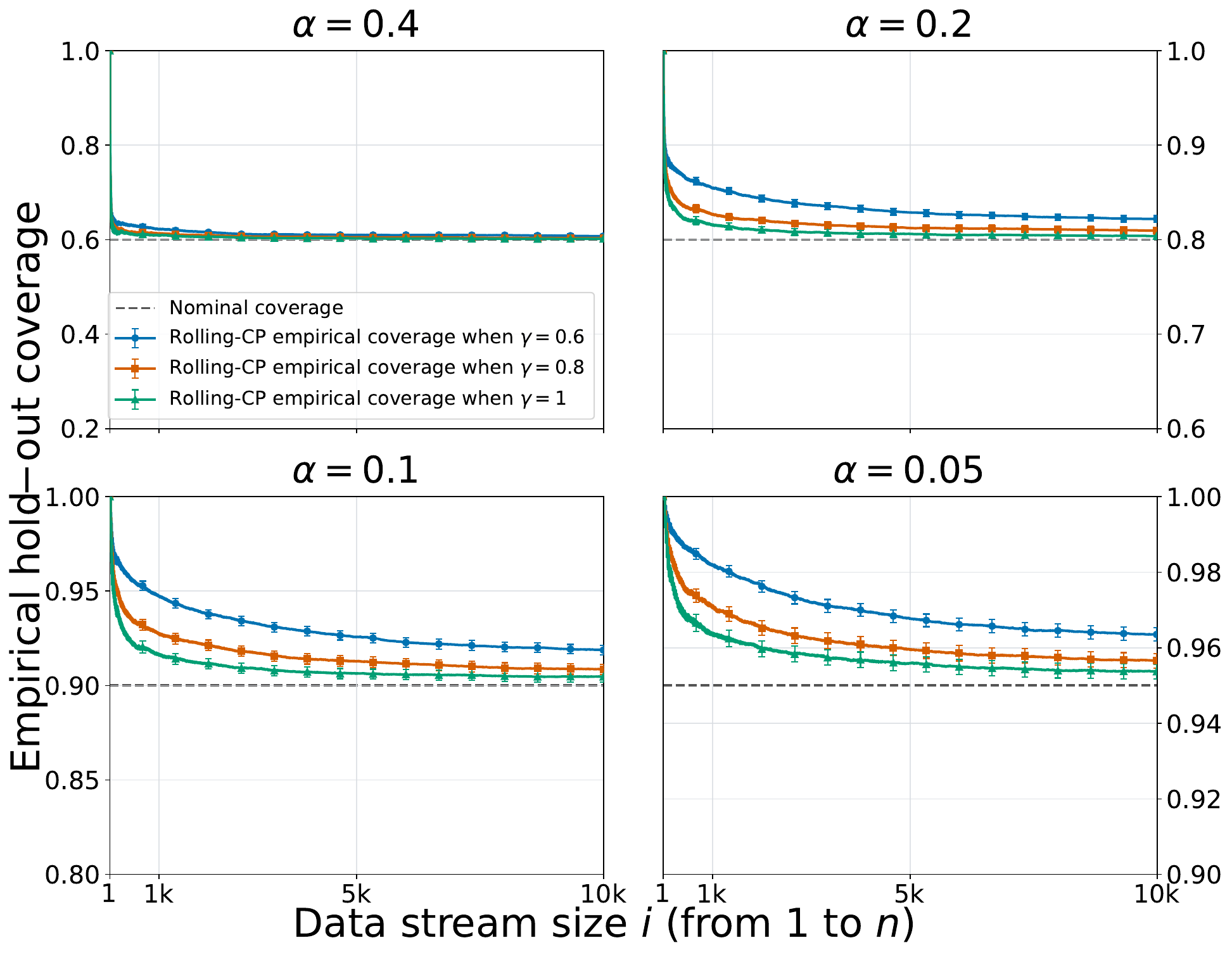}
\par\vspace{0.1em}
\parbox{\linewidth}{\footnotesize \textbf{(c)} Average held-out coverage along the data stream for $s_i^{\mathsf{Ent}}$ and $\alpha\in\{0.4,0.2,0.1,0.05\}$.}
\end{minipage}%
\hspace{0.03\textwidth}
\begin{minipage}[t]{0.45\textwidth}
\centering
\includegraphics[width=1.02\linewidth]{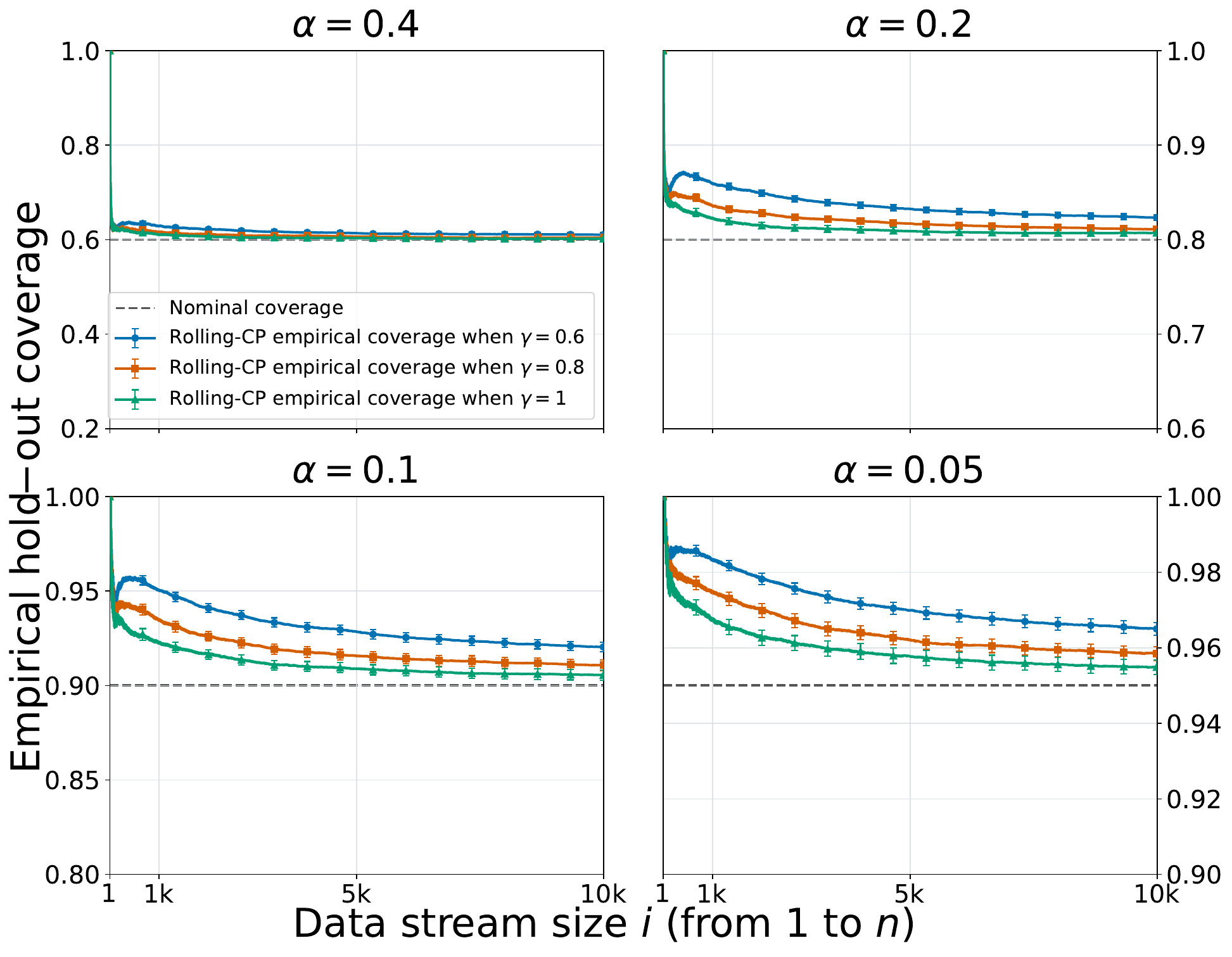}
\par\vspace{0.1em}
\parbox{\linewidth}{\footnotesize \textbf{(d)} Average held-out coverage along the data stream for $s_i^{\mathsf{Mar}}$ and $\alpha\in\{0.4,0.2,0.1,0.05\}$.}
\end{minipage}
\par\vspace{0.5em}

\caption{Numerical experiments of SGD for multi-index logistic regression.
\textbf{(a)} End-of-stream empirical marginal coverage using $s_i^{\mathsf{Ent}}$.
\textbf{(b)} End-of-stream empirical marginal coverage using $s_i^{\mathsf{Mar}}$.
\textbf{(c)} Evolution of coverage using $s_i^{\mathsf{Ent}}$.
\textbf{(d)} Evolution of coverage using $s_i^{\mathsf{Mar}}$.
All panels use $d=10$, $n=10,000$, $M=100$, $\eta_0=1$, $t_0=10$, $T=100$, and $\gamma\in\{0.6,0.8,1\}$.
Error bars are pointwise $95\%$ confidence intervals across trials.
Panels \textbf{(e)} and \textbf{(f)} appear in the continuation of this figure.}
\label{fig:SGD-logistic}
\end{figure}

\begin{figure}[!htbp]
\ContinuedFloat
\centering
\begin{minipage}[t]{0.45\textwidth}
\centering
\includegraphics[width=1.02\linewidth]{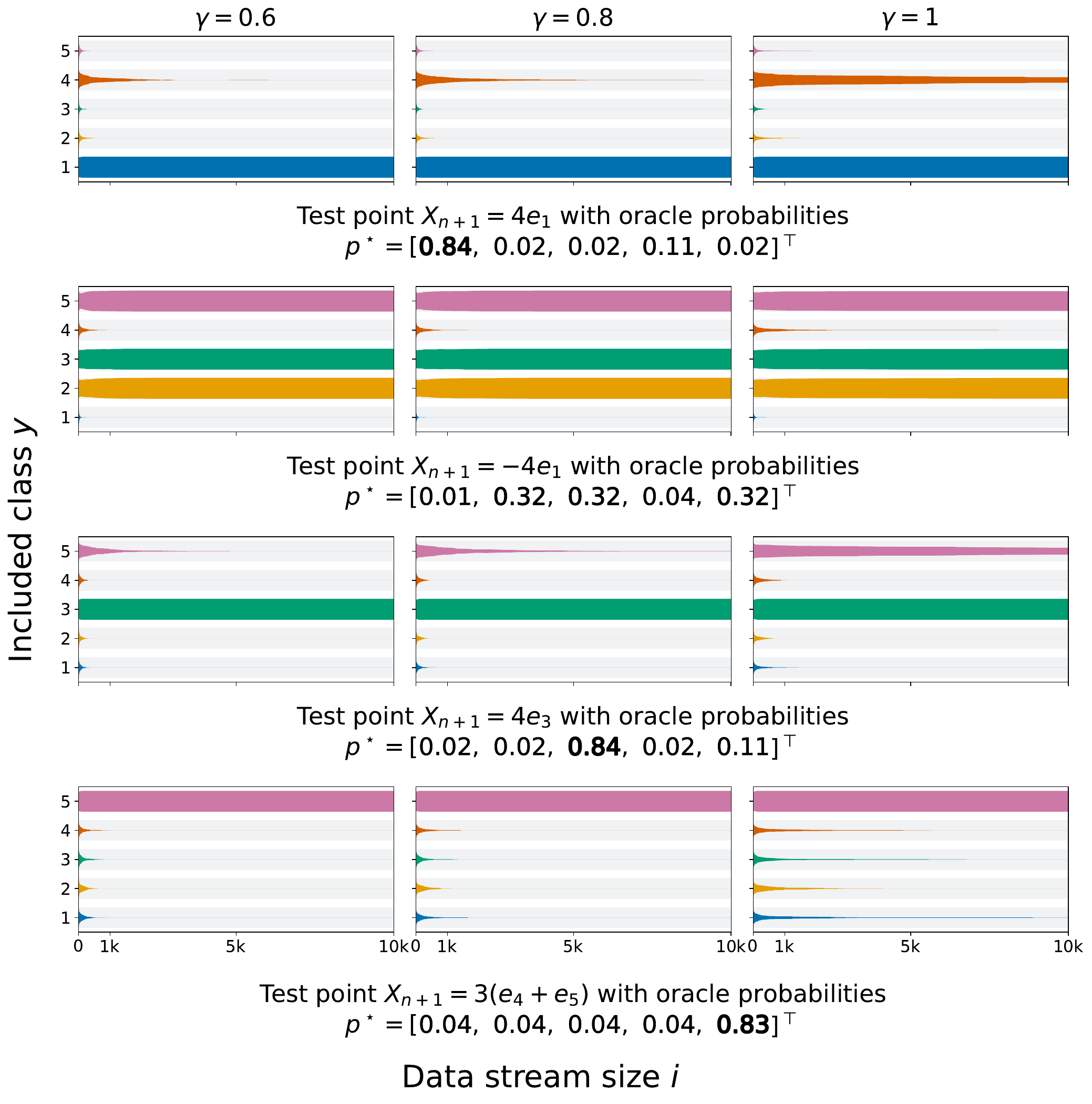}
\par\vspace{0.1em}
\parbox{\linewidth}{\footnotesize \textbf{(e)} Class-wise inclusion frequencies of the $70\%$ rolling-CP set for four fixed test features using $s_i^{\mathsf{Ent}}$.}
\end{minipage}%
\hspace{0.03\textwidth}
\begin{minipage}[t]{0.45\textwidth}
\centering
\includegraphics[width=1.02\linewidth]{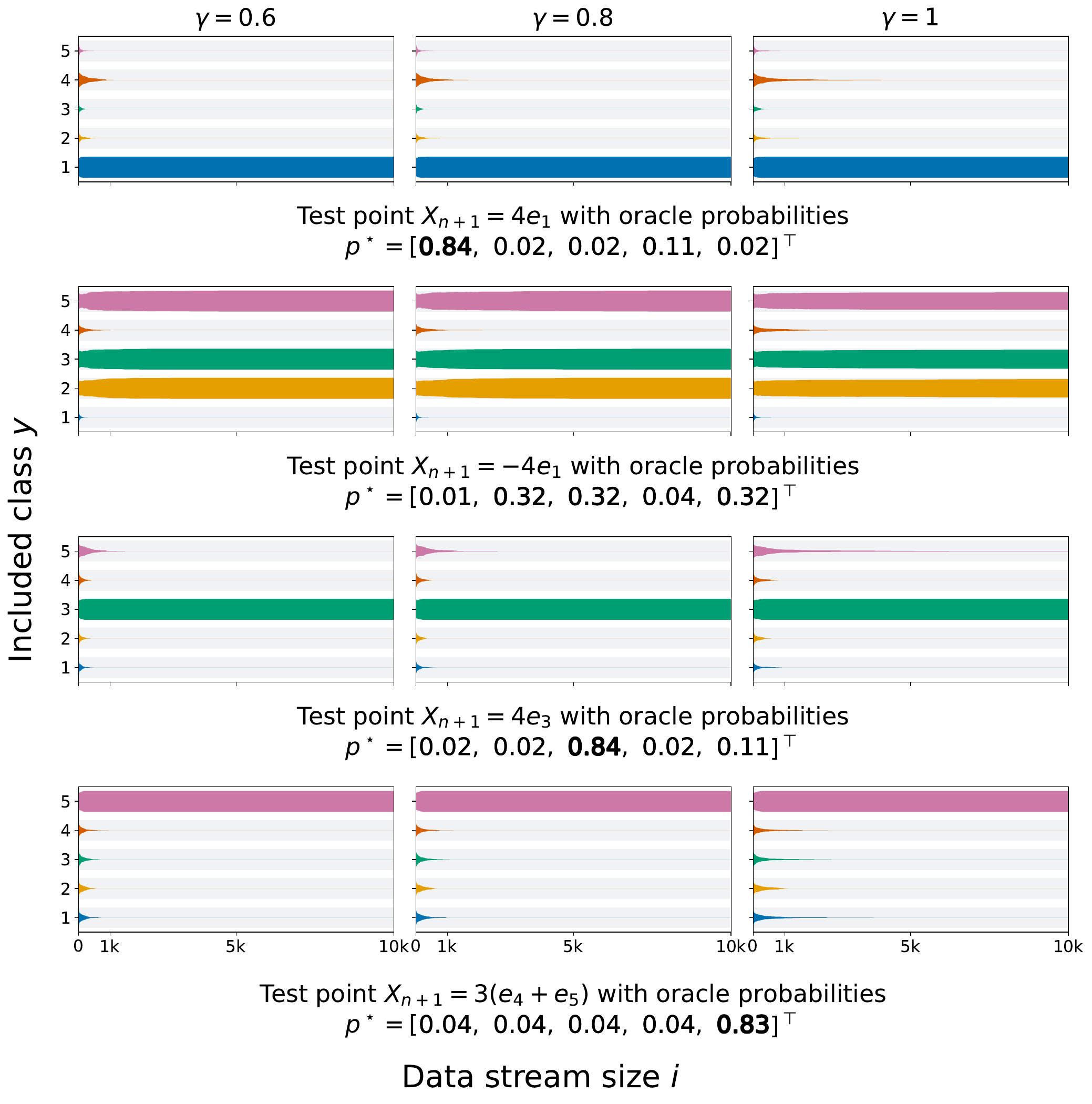}
\par\vspace{0.1em}
\parbox{\linewidth}{\footnotesize \textbf{(f)} Class-wise inclusion frequencies of the $70\%$ rolling-CP set for four fixed test features using $s_i^{\mathsf{Mar}}$.}
\end{minipage}
\par\vspace{0.5em}

\caption[]{Numerical experiments of SGD for multi-index logistic regression (continued).
\textbf{(e)} Fixed-feature class-inclusion paths using $s_i^{\mathsf{Ent}}$.
\textbf{(f)} Fixed-feature class-inclusion paths using $s_i^{\mathsf{Mar}}$.
Both panels use $d=10$, $n=10,000$, $M=100$, $\eta_0=1$, $t_0=10$, $T=100$, and $\gamma\in\{0.6,0.8,1\}$.
We examine the output of the prediction set $\mc{C}_n(x)$ at a fixed test point feature $x$, as $n$ increases; boldface in the test-point annotations reflects the magnitudes of the oracle class probabilities for this value of $x$.}
\end{figure}

\subsection{One-pass training on MNIST}\label{sec:numerical-mnist}

We finally consider the MNIST handwritten-digit benchmark, consisting of $60,000$ training images and $10,000$ test images from ten classes~\citep{lecun1998gradient}. We randomly order the training observations and process each observation exactly once. Our predictor is a shallow four-layer LeNet-style convolutional network. We refer readers to \citet{lecun1998gradient} for background on the LeNet architecture. The precise network architecture and complete experimental setup are detailed in the \href{https://github.com/Moriartycc/rolling-conformal/tree/main/one-pass-mnist}{code repository}.

Let $f_{\Theta}(x)\in\R^{K}$ with $K=10$ denote the network logits, where $\Theta$ collects all trainable parameters of the neural network. Starting from $\what{\Theta}_0$, we train by pure SGD with batch size one:
\begin{align*}
    \what{\Theta}_i = \what{\Theta}_{i-1} - \eta_i \left.
    \nabla_{\Theta}
    \brc{-\log\mathsf{softmax}\prn{f_{\Theta}(X_i)}_{Y_i}}\right|_{\Theta=\what{\Theta}_{i-1}}, \qquad
    \eta_i = \frac{\eta_0}{(t_0+i)^\gamma}.
\end{align*}
We take $\eta_0=50, t_0=5,000$ and $\gamma = 1$ for this experiment. The gradients are computed using PyTorch's built-in automatic differentiation and backpropagation, and the complete parameter vector is updated according to the displayed SGD recursion. For $z=(x,y)$, we construct rolling-CP using the cross-entropy score
\begin{align*}
    s_i(z;Z_{<i})
    =
    -\log\mathsf{softmax}\prn{f_{\what{\Theta}_{i-1}}(x)}_y,
\end{align*}
evaluated immediately before the SGD update at iteration $i$.

Unlike the preceding experiments, we run a single training trajectory and estimate coverage using $M=1,000$ fixed observations drawn from the MNIST test set. Thus, the reported inclusion frequencies average over held-out test observations conditional on the realized training stream and algorithmic initialization. Consequently, this experiment does not fully evaluate the marginal guarantees in Theorems~\ref{thm:exchangeable-marginal} and~\ref{thm:cvx-dcx}, which also average over the randomness of the training data and algorithm. A direct numerical verification of these marginal guarantees would require multiple independent training streams of comparable size, which is infeasible given that the single standard MNIST training set contains only $60,000$ observations. Further dividing it would no longer represent the intended large-scale, one-pass training regime. We therefore retain the complete training stream and use the held-out test observations to estimate coverage conditional on this realized trajectory. 

\begin{figure}[!htbp]
    \centering
    \begin{minipage}[t]{0.55\textwidth}
        \centering
        \includegraphics[width=0.95\linewidth]{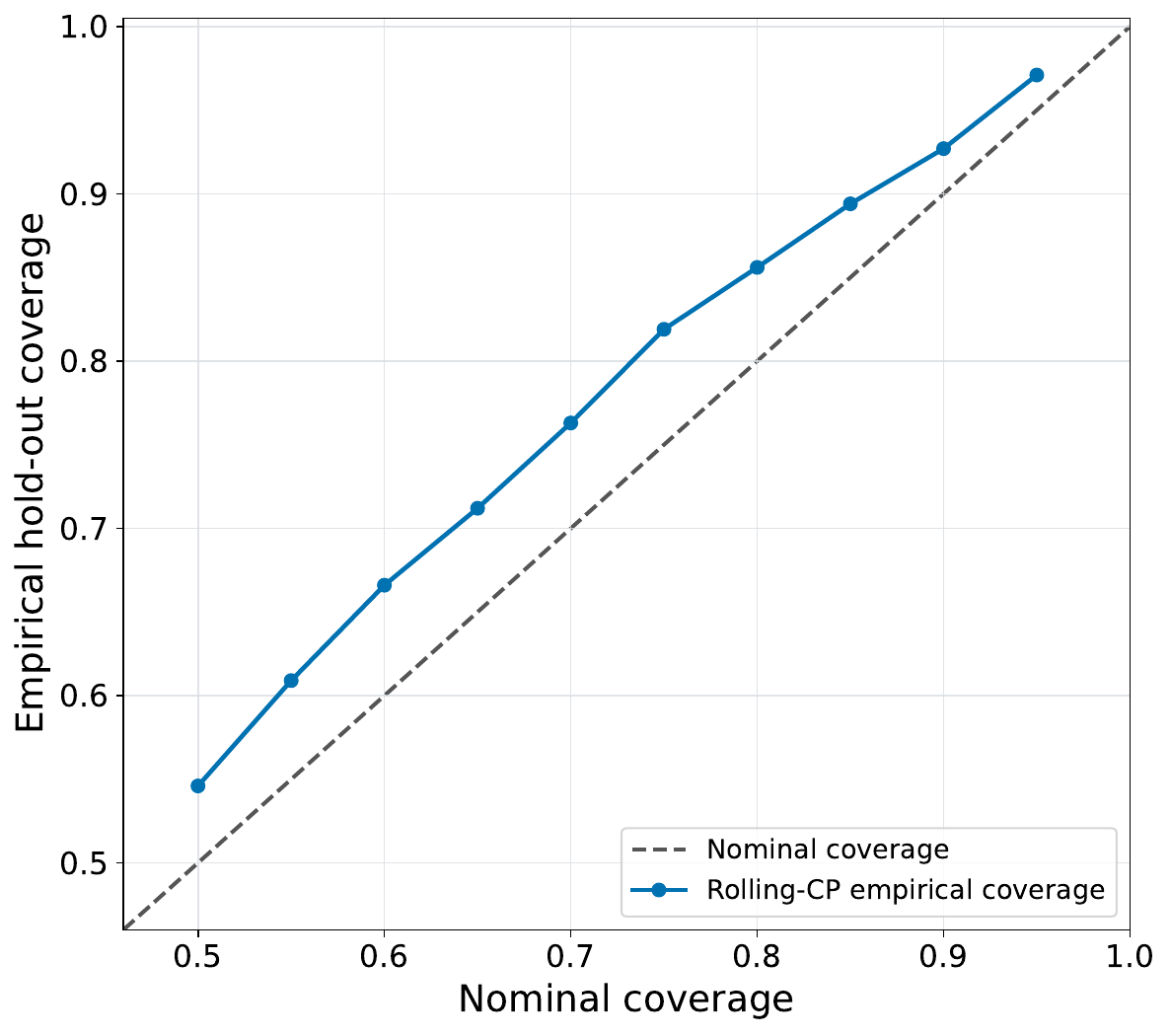}
        \par\vspace{0.2em}
        \parbox{\linewidth}{\footnotesize
        \textbf{(a)} For each nominal level, the end-of-stream coverage is the inclusion frequency among the $M=1,000$ held-out test observations.}
    \end{minipage}%
    
    \centering
    \begin{minipage}[t]{0.55\textwidth}
        \centering
        \includegraphics[width=0.95\linewidth]{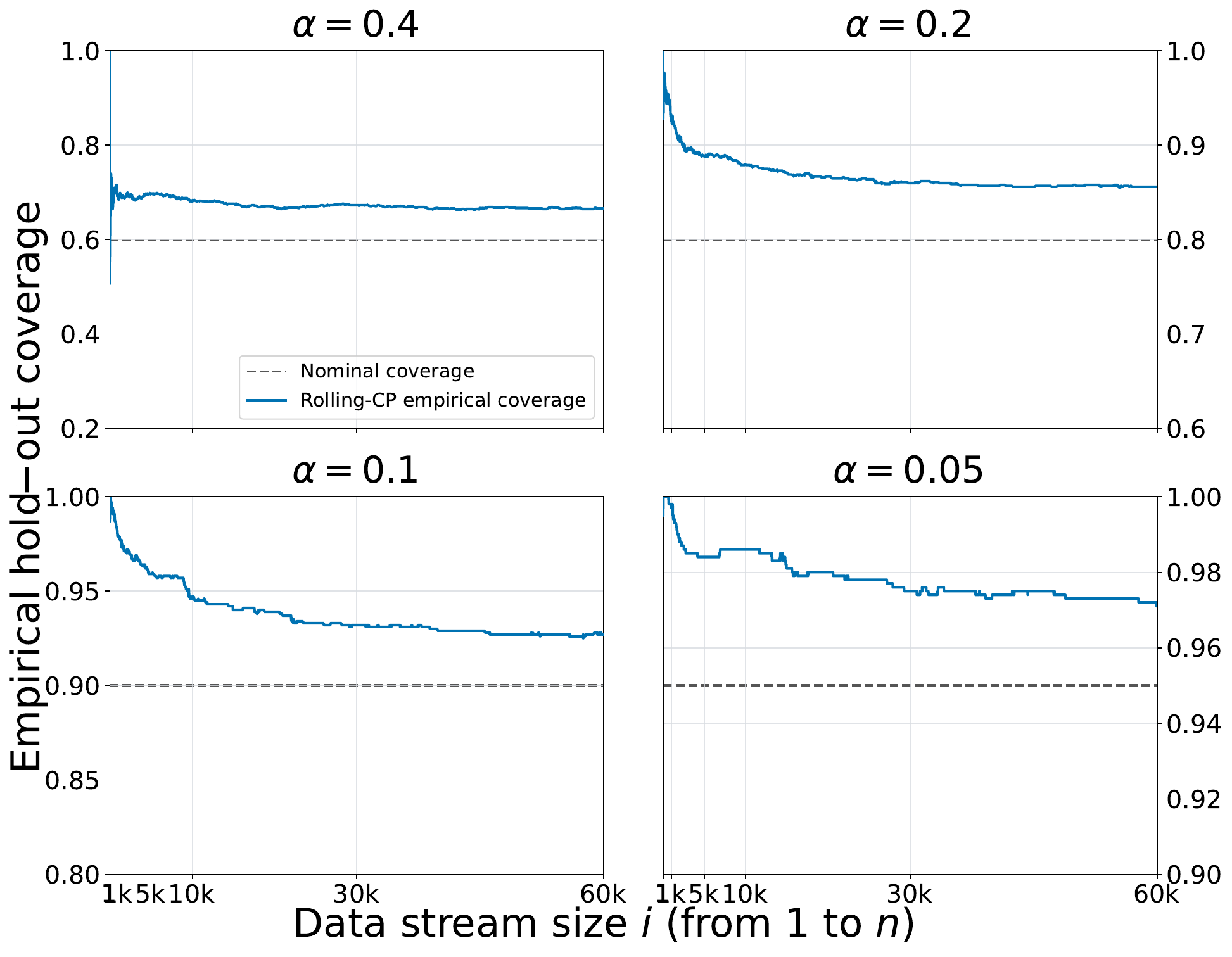}
        \par\vspace{0.2em}
        \parbox{\linewidth}{\footnotesize
        \textbf{(b)} For $\alpha\in\{0.4,0.2,0.1,0.05\}$, we report the held-out inclusion frequency along the training stream.}
    \end{minipage}

    \par\vspace{0.5em}
    \caption{Numerical experiments of one-pass SGD on MNIST. \textbf{(a)} End-of-stream empirical held-out coverage versus nominal coverage. \textbf{(b)} Evolution of empirical held-out coverage during the single training pass. Both panels use $n=60,000$ training observations, $M=1,000$ held-out test observations, which were drawn independently from the test set.}
    \label{fig:one-pass-MNIST}
\end{figure}

Consequently, the experiment is more closely related to the quantity of training-conditional coverage controlled in Theorem~\ref{thm:training-conditional}, although a full numerical verification of its high-probability guarantee would still require repetition over independent training streams. Despite these limitations, Figure~\ref{fig:one-pass-MNIST} already exhibits reasonable training-conditional coverage along the realized training trajectory.

%% file: discussion.tex
\section{Discussion}\label{sec:discussion}

In this work, we introduce rolling-CP, a novel distribution-free predictive inference method tailored to sequential training settings that play a central role in modern machine learning. Rolling-CP is computationally lean and tuning-free, and naturally adapts to evolving models. It provides validity under minimal distributional assumptions, without requiring additional data splitting or repeated retraining. 

Next, we discuss several interesting open questions and directions for future work. A promising direction is to use rolling-CP to curate training and track qualitative changes along the trajectory. Classical early stopping can be statistically optimal in some regimes \citep{raskutti2014early}, while modern overparameterized models can show double descent \citep{nakkiran2020deep} or grokking \citep{power2022grokking} (where generalization improves only long after the training data are fit). These arise in different regimes and imply there is no universal stopping rule based solely on training loss or iteration count. Because rolling-CP records each example's score before using that example in the next update, it supplies an out-of-sample diagnostic within a single run. These diagnostics can guide stopping or continuing training, and help distinguish overfitting from delayed generalization and other training transitions.

It is also interesting to extend rolling-CP from marginal and training-conditional validity to test-conditional coverage, ensuring reliability at a given feature value rather than only on average over the test population. Exact, nontrivial, distribution-free test-conditional coverage is generally impossible without additional assumptions \citep{vovk2012conditional,barber2021limits,areces2024fundamental}. However, it is still possible to relax to obtain coverage over structured groups or weighting functions \citep{gibbs2025conditional} and randomized localization around the test point \citep{hore2025localweights}. The technical challenge remains extending these ideas to rolling-CP requires simultaneously accounting for scores that evolve throughout training.

Another potential direction is to study the performance of rolling-CP under temporally dependent data streams, where exchangeability no longer holds. It would be interesting to understand whether rolling-CP remain empirically robust in such time series setting, as split conformal prediction often does under short-range temporal dependence \citep{barber2025predictive}, and admits a corresponding theoretical guarantee, or it can suffer severe loss of coverage and requires further modification, as in related extensions of jackknife+ \citep{jiang2026leave}. We leave these questions for future work.

%% file: proof-main-results.tex
\section{Proofs for Section~\ref{sec:main-results}} \label{proof:main-results}

\subsection{Proof of Theorem~\ref{thm:cvx-dcx}} \label{proof:cvx-dcx}
As discussed in Section~\ref{sec:cvx}, the proof of our main coverage guarantee (Theorem~\ref{thm:exchangeable-marginal}) relies on the convex ordering result stated in Theorem~\ref{thm:cvx-dcx}. In this section, we complete the picture by providing the proof of Theorem~\ref{thm:cvx-dcx}.
Before presenting the proof of the theorem, we first introduce some additional background and definitions.

In the proof, we will also use conditional versions of the orders $\stl, \dcx, \cvx$. For example,
\begin{align*}
(X\mid Z)\dcx(Y\mid Z)
\end{align*}
means that, for every nonincreasing convex function $f$ for which the conditional expectations are well-defined,
\begin{align*}
\E[f(X)\mid Z]\leq \E[f(Y)\mid Z]\qquad\textnormal{almost surely}.
\end{align*}
The conditional stochastic and convex orders are defined analogously. A useful equivalent representation of the decreasing convex order is the following \citep{wang2024testing}:
\begin{equation}\label{eqn:dcx_equiv}
X\dcx Y \quad \Longleftrightarrow \quad \begin{array}{c}\textnormal{ there exist jointly distributed random variables $(X',Y')$,}\\\textnormal{with $X'\eqd X$ and $Y'\eqd Y$, such that $\mathbb{E}[Y'\mid X'] \leq X'$ almost surely.}\end{array}
\end{equation}
Following the terminology and results of \citet{wang2024testing}, a random variable $p\in[0,1]$ satisfying\footnote{Throughout, we will use the notation of random variables or of their distributions, as convenient. Specifically, the statement $p\dcx \mathsf{Unif}(0,1)$ should be interpreted as equivalent to stating $p\dcx U$ where $U\sim\mathsf{Unif}(0,1)$.} $p\dcx \mathsf{Unif}(0,1)$ is called a ``p*-value'', and satisfies the following guarantee:
\begin{equation}\label{eqn:factor_of_2__p*_variable}
\textnormal{If $p\in[0,1]$ satisfies $p\dcx\mathsf{Unif}(0,1)$, then }\P \prn{p\leq \alpha}\leq 2\alpha\textnormal{ for all $\alpha\in[0,1]$.}
\end{equation}

Next we need another definition:
\begin{definition}[Comonotone random variables]
    Jointly distributed random variables $X,Y\in\R$ are \emph{comonotone} if there exists a random variable $Z$ and nondecreasing functions $f,g:\R\to\R$ such that $X=f(Z)$ and $Y=g(Z)$, almost surely. Moreover, for any random variables $X,Y$, there exists a unique \emph{comonotone coupling}: a joint distribution on $(X^*,Y^*)$, such that $X^*,Y^*$ are comonotone, and such that marginally it holds that $X^*\eqd X$ and $Y^*\eqd Y$. 
\end{definition}

We will need one additional result that relies on the decreasing convex order.
\begin{lemma}\label{lem:dcx_comonotone}
    Let $(A,B)$ be jointly distributed random variables, and let $X,Y$ be random variables such that $A\dcx X$ and $B\dcx Y$. Let $(X^*,Y^*)$ be a comonotone coupling of $X$ and $Y$. Then
    \[A+B\dcx X^*+Y^*.\]
\end{lemma}

We are now ready to prove the theorem.
\paragraph{Proof of Theorem~\ref{thm:cvx-dcx}.}
We prove the result by induction on $n$. 
For $n=0$, it holds trivially as $p_{\mathsf{rolling}}=1$ almost surely. Now assume $n\geq 1$ and suppose the theorem is true with $n-1$ in place of $n$. We will now condition on $Z_1$ and on the empirical distribution $\what{P}_Z = \frac{1}{n+1}\sum_{i=1}^{n+1}\delta_{Z_i}$ of the sequence. Let $\wt{Z}_i = Z_{i+1}$ for $i=1,\dots, n$, and define functions
\[\wt{s}_i(z;(z_1,\dots,z_{i-1})) = s_{i+1}(z;(Z_1,z_1,\dots,z_{i-1})),\]
for $i=1,\dots,n-1$ (note that these functions depend implicitly on $Z_1$; since we will condition on $Z_1$, we do not make this argument explicit).
Define similarly $\wt{Z}_{<i} = (\wt{Z}_1, \dots, \wt{Z}_{i-1})$ and 
\[\wt{p}_{\mathsf{rolling}} = \frac{1 + \sum_{i=2}^n  \ind\{s_i(Z_i; Z_{<i}) \geq s_i(Z_{n+1}; Z_{<i})\}}{n} = \frac{1 + \sum_{i=1}^{n-1}  \ind\{\wt{s}_i(\wt{Z}_i; \wt{Z}_{<i}) \geq  \wt{s}_i(\wt{Z}_n;\wt{Z}_{<i})\}}{n}.\]
By exchangeability of $(Z_1,Z_2\dots,Z_{n+1}) = (Z_1,\wt{Z}_1,\dots,\wt{Z}_n)$, it holds that $\wt{Z}_1,\dots,\wt{Z}_n$ are exchangeable conditional on $(Z_1,\what{P}_Z)$. Therefore, applying the induction hypothesis (with $n-1$ in place of $n$, and with $\wt{Z}_i,\wt{s}_i$ in place of $Z_i,s_i$), we have
\[(\wt{p}_{\mathsf{rolling}}\mid Z_1,\what{P}_Z) \dcx \mathsf{Unif}\left(\left\{\frac{1}{n},\dots,\frac{n-1}{n},1\right\}\right).\]
Now write
\[p_{\mathsf{rolling}} = \frac{ 1 + \ind\{s_1(Z_1) \geq s_1(Z_{n+1})\} + \sum_{i=2}^n  \ind\{s_i(Z_i;Z_{<i})\geq s_i(Z_{n+1};Z_{<i})\}}{n+1} =  \frac{A + n\cdot \wt{p}_{\mathsf{rolling}}}{n+1}\]
where 
\[A = \ind\{s_1(Z_1) \geq s_1(Z_{n+1})\}.\]
(Note that since $Z_{<1}$ contains no information, we write $s_1(z)$ in place of $s_1(z; Z_{<1})$.)
Define $h(Z_1,\what{P}_Z)\in[0,1]$ as 
\[h(Z_1,\what{P}_Z):= 1 - \mathbb{E}[A \mid Z_1,\what{P}_Z] = \mathbb{P}\prn{s_1(Z_1)<s_1(Z_{n+1})\mid Z_1,\what{P}_Z} = \frac{1}{n}\sum_{i=2}^{n+1}\ind\{s_1(Z_1)<s_1(Z_i)\},\] 
where the last step holds since $Z_2,\dots,Z_{n+1}$ are exchangeable conditionally on $Z_1,\what{P}_Z$. Then,
conditional on $(Z_1,\what{P}_Z)$, we have
\[A\mid Z_1,\what{P}_Z \sim \textnormal{Bernoulli}(1-h(Z_1,\what{P}_Z)), \quad n\cdot \wt{p}_{\mathsf{rolling}}\mid (Z_1,\what{P}_Z) \dcx B\]
where we define the random variable $B\sim \mathsf{Unif}([n])$. Next we construct a comonotone coupling of the conditional distributions of $A\mid Z_1,\what{P}_Z$ and $B\mid Z_1,\what{P}_Z$: define
\[A^* = \ind\{U>  h(Z_1,\what{P}_Z)\}, \quad B^* = \lceil n \cdot U \rceil,\]
where $U\sim\mathsf{Unif}(0,1)$. Applying Lemma~\ref{lem:dcx_comonotone}, we have
\[\left(A + n\cdot \wt{p}_{\mathsf{rolling}}\,\middle| \, Z_1,\what{P}_Z \right)\dcx \left( A^*+B^* \,\middle|\, Z_1,\what{P}_Z\right).\]
After rescaling, and plugging in the construction of $(A^*,B^*)$, this yields
\[(p_{\mathsf{rolling}}\mid Z_1, \what{P}_Z) \dcx\left( \frac{ \ind\{U>h(Z_1,\what{P}_Z)\} + \lceil n \cdot U \rceil}{n+1}\,\middle|\, Z_1, \what{P}_Z\right).\]
Since the decreasing convex order is preserved under marginalization, we have therefore proved that
\[(p_{\mathsf{rolling}} \mid \what{P}_Z) \dcx \left( \frac{ \ind\{U>h(Z_1,\what{P}_Z)\} + \lceil n \cdot U \rceil}{n+1} \,\middle|\, \what{P}_Z\right),\]
where we have $U\sim\mathsf{Unif}(0,1)$ drawn independently of $(Z_1,\what{P}_Z)$. Next, by exchangeability of $Z_1,\dots,Z_{n+1}$, we can verify that
\[(h(Z_1,\what{P}_Z) \mid \what{P}_Z)\stl C\textnormal{ where }C\sim \mathsf{Unif}\left(\left\{0,\frac{1}{n},\dots,\frac{n-1}{n},1\right\}\right).\]
Since $$\frac{ \ind\{U > h(Z_1,\what{P}_Z)\} + \lceil n \cdot U \rceil}{n+1}$$ is nonincreasing as a function of $h(Z_1,\what{P}_Z)$, and $U$ is independent from $Z_1,\what{P}_Z$,  this means that
\[\frac{ \ind\{U> h(Z_1,\what{P}_Z)\} + \lceil n \cdot U \rceil}{n+1} \stg \frac{ \ind\{U> C\} + \lceil n \cdot U \rceil}{n+1},\]
where $U$ and $C$ are independent.
Thus, since $X\succeq_{\mathrm{st}}Y$ implies $X\dcx Y$ we have
\[\frac{ \ind\{U> h(Z_1,\what{P}_Z)\} + \lceil n \cdot U \rceil}{n+1} \dcx \frac{ \ind\{U> C\} + \lceil n \cdot U \rceil}{n+1},\]
Combined with the work above, therefore,
\[p_{\mathsf{rolling}} \dcx \frac{ \ind\{U> C\} + \lceil n \cdot U \rceil}{n+1}.\]Finally, a straightforward calculation shows that $\ind\{U> C\} + \lceil n \cdot U \rceil\sim\mathsf{Unif}([n+1])$, which completes the proof.
\endproof

\subsection{Alternative proof of Theorem~\ref{thm:exchangeable-marginal}} \label{proof:exchangeable-marginal}

Earlier, in Section~\ref{sec:cvx}, we proved Theorem~\ref{thm:exchangeable-marginal} in the general case of exchangeable data by first establishing a result on distributional ordering for the rolling-CP set. Here, we give an alternative proof that is more direct, without relying on distributional ordering; this proof is more similar to the original intuitive argument given for the special case of i.i.d.\ data, in Section~\ref{sec:marginal-coverage}.

We will prove the result again assuming that there are no ties among scores, almost surely. (Note that the proof given in Section~\ref{sec:cvx} did not require this assumption; here we are only presenting an alternative argument for intuition and thus will make this simplifying assumption.) We will implicitly condition on the empirical distribution $\what{P}_Z = \frac{1}{n+1}\sum_{i=1}^{n+1} \delta_{Z_i}$ throughout.

We use the fact conditional on $Z_{<i}$, the future data stream $Z_i, \dots, Z_n, Z_{n+1}$ are exchangeable. For $i=1,2, \dots, n$, define
\begin{align*}
    p_i = \frac{\sum_{j=i+1}^{n+1} \ind \brc{s_i(Z_j; Z_{<i}) > s_i(Z_i; Z_{<i})}}{n-i+1} =: \frac{r_i-1}{n-i+1},
\end{align*}
Here, $r_i = r_i(Z_i; Z_{<i})$ denotes the conditional order statistic of $Z_i$ within the collection $(Z_i, \dots, Z_{n+1})$. 
We note that there is no ``$+1$'' term in this construction, and hence this is not a valid p-value---indeed, we may have $p_i=0$. By definition, this quantity follows the distribution
\begin{align*}
    p_i \sim \mathsf{Unif}\prn{\brc{0, \frac{1}{n-i+1},\dots,\frac{n-i}{n-i+1}, 1}}.
\end{align*}

We next claim that $p_i$'s are indeed mutually independent (but importantly, not identically distributed). To see that, we first notice that by exchangeability of $(Z_i, \dots, Z_n, Z_{n+1})$ and that $r_1,\dots, r_{i-1}$ are all invariant under any permutation of $(Z_i, \dots, Z_n, Z_{n+1})$,
\begin{align*}
    r_i \mid Z_{<i}, r_1, \dots, r_{i-1} \sim \mathsf{Unif}([n-i+2]),
\end{align*}
which implies $r_1,\dots,r_n$ are mutually independent, by induction; this is sufficient since $p_i$ is a function of $r_i$, for each $i$.

Next, define events
\[\mc{A}_i = \brc{s_i(Z_{n+1};Z_{<i})>s_i(Z_i;Z_{<i})}.\]
By exchangeability of the data, it holds that
\[\P\prn{\mc{A}_i\mid Z_{<i+1}} = p_i,\]
since conditional on $Z_{<i+1} = (Z_1,\dots,Z_i)$, the remaining data points $Z_{i+1},\dots,Z_{n+1}$ are exchangeable. 

Now we need a lemma:
\begin{lemma}\label{lem:pigeonhole_exch}
    Let $\mc{A}_1,\dots,\mc{A}_n$ be any events. Let $\mc{F}_1\subseteq \dots\subseteq \mc{F}_n$ be some filtration, and define the random variables 
    \[p_i = \P\prn{\mc{A}_i\mid \mc{F}_i}\in\mc{F}_i.\]
    Suppose $p_1,\dots,p_n$ are mutually independent, and have marginal distributions
    \[p_i\sim\mathsf{Unif}\left(\left\{0, \frac{1}{n-i+1},\dots,\frac{n-i}{n-i+1},1\right\}\right).\]
    Then for any $\Delta \geq 0$,
    \[\P\prn{\sum_{i=1}^n \ind_{\mc{A}_i} \geq  n-\Delta} \leq \frac{2\Delta +1}{n+1}.\]
\end{lemma}
Applying the lemma with $\Delta = \alpha(n+1)-1$, and $\mc{F}_i$ defined as the $\sigma$-algebra generated by $Z_{<i+1} = (Z_1,\dots,Z_i)$, we therefore have
\[\P\prn{\sum_{i=1}^n \ind_{\mc{A}_i}\geq (1-\alpha)(n+1)}\leq 2\alpha.\]
As before, by definition of the rolling-CP prediction set, this is sufficient, since
\[Z_{n+1}\not\in\mc{C}_n \ \Longleftrightarrow \ \sum_{i=1}^n\ind_{\mc{A}_i} \geq (1-\alpha)(n+1).\]

\subsection{Proof of Theorem~\ref{thm:training-conditional}}
     Let $U_1,\cdots, U_n \simiid \mathsf{Unif}(0,1)$ be drawn independently of the data stream $Z_1,\cdots, Z_n, Z_{n+1}$. For each $i \in [n]$, define
    \[p_i = U_i \cdot \P(s_i(Z_{n+1};Z_{<i}) = s_i(Z_i;Z_{<i})\mid Z_{<i+1}) +   \P(s_i(Z_{n+1};Z_{<i}) > s_i(Z_i;Z_{<i})\mid Z_{<i+1}). \]
    Note that when the score comparisons are almost surely free of ties, $p_i$ defined above essentially is the same as it was defined in the proof of marginal coverage in the i.i.d.\ case in Section~\ref{sec:marginal-coverage}. Hence, following the same argument as in Section~\ref{sec:marginal-coverage}, it follows inductively that
    \begin{align*}
    p_1,\dots,p_n\simiid\mathsf{Unif}(0,1).
    \end{align*}
    Note also that $(p_1,\dots,p_n)$ is a function of $Z_{<n+1} = (Z_1,\dots,Z_n)$, and is therefore independent of $Z_{n+1}$.

Next, define events
\[\mc{A}_i = \brc{s_i(Z_{n+1};Z_{<i})>s_i(Z_i;Z_{<i})}.\]
Then
\[\P\prn{\mc{A}_i\mid Z_{<n+1}} =\P\prn{\mc{A}_i\mid Z_{<i+1}} \leq p_i, \]
where the first step holds since $\mc{A}_i$ depends only $Z_1,\dots,Z_i,Z_{n+1}$ and is independent of $Z_{i+1},\dots,Z_n$ (exactly as in the proof of Theorem~\ref{thm:exchangeable-marginal} for the special case of i.i.d.\ data, shown in Section~\ref{sec:marginal-coverage}).

    Our next step is a lemma, similar to the result of Lemma~\ref{lem:pigeonhole}:
    \begin{lemma}\label{lem:pigeonhole_conditional}
    Let $\mc{A}_1,\dots,\mc{A}_n$ be any events. Suppose that there exists some random variable $W$, and random variables $p_1,\dots,p_n$, such that 
    \[\P\prn{\mc{A}_i \mid W}\leq p_i\textnormal{ almost surely for all $i\in[n]$, and }p_1,\dots,p_n\simiid \mathsf{Unif}(0,1).\]
    Then for any $\Delta,\delta \geq 0$,
    \[\P\prn{\sum_{i=1}^n \ind_{\mc{A}_i} \geq  n-\Delta \,\bigg|\, W} \leq \frac{2\Delta+1}{n+1} +  \sqrt{\frac{n\log(1/\delta)}{2(\lfloor\Delta\rfloor+1)^2}} \textnormal{ with probability $\geq 1-\delta$}.\]
\end{lemma}
Since $Z_{n+1}\in\mc{C}_n$ if and only if $\sum_{i=1}^n \ind_{\mc{A}_i}<(1-\alpha)(n+1)$, we therefore have
\[\alpha_P(Z_{<n+1}) = \P\prn{\sum_{i=1}^n \ind_{\mc{A}_i} \geq (1-\alpha)(n+1)\,\bigg|\, Z_{<n+1}} \leq 2\alpha + \sqrt{\frac{n\log(1/\delta)}{2\lfloor\alpha(n+1)\rfloor^2}},\]
with probability $\geq 1-\delta$, by applying the above lemma with $\Delta = \alpha(n+1)-1$ and $W=Z_{<n+1}$. Finally, if $\alpha(n+1)\geq 1$, then we have $\sqrt{\frac{n\log(1/\delta)}{2\lfloor\alpha(n+1)\rfloor^2}}\leq \sqrt{\frac{n\log(1/\delta)}{2\left(\frac{\alpha(n+1)}{2} \right)^2}} \leq \sqrt{\frac{2\log(1/\delta)}{\alpha^2(n+1)}}$, while if instead $\alpha(n+1)<1$ then we must have $\alpha_P(Z_{<n+1})=0$ almost surely (since $p_{\mathsf{rolling}}\geq \frac{1}{n+1}$ almost surely, by construction), which proves the desired bound.

    Finally, to obtain the uniform bound, it suffices to note that for any $\delta \in (0,1)$, $\sum_{n \ge 1} (\delta/(n+1)^2) \le \delta$. By the work above (with $\delta/(n+1)^2$ in place of $\delta$), for each $n\geq 1$ the event
    \[\alpha_P(Z_{< n+1}) \le 2\alpha +\sqrt{\frac{2\log((n+1)^2/\delta)}{\alpha^2(n+1)}}\]
    holds with probability $\geq 1-\delta/(n+1)^2$. Therefore, taking a union bound over these events for $n\geq 1$ (along with the fact that $\alpha_P(Z_{<1}) = 0$ almost surely, since $p_{\mathsf{rolling}}=1$ almost surely at time $n=0$), yields the desired uniform bound over all $n\geq 0$.

    \subsection{Proof of Proposition~\ref{prop:tightness}} \label{proof:tightness}
\begin{figure}[t]
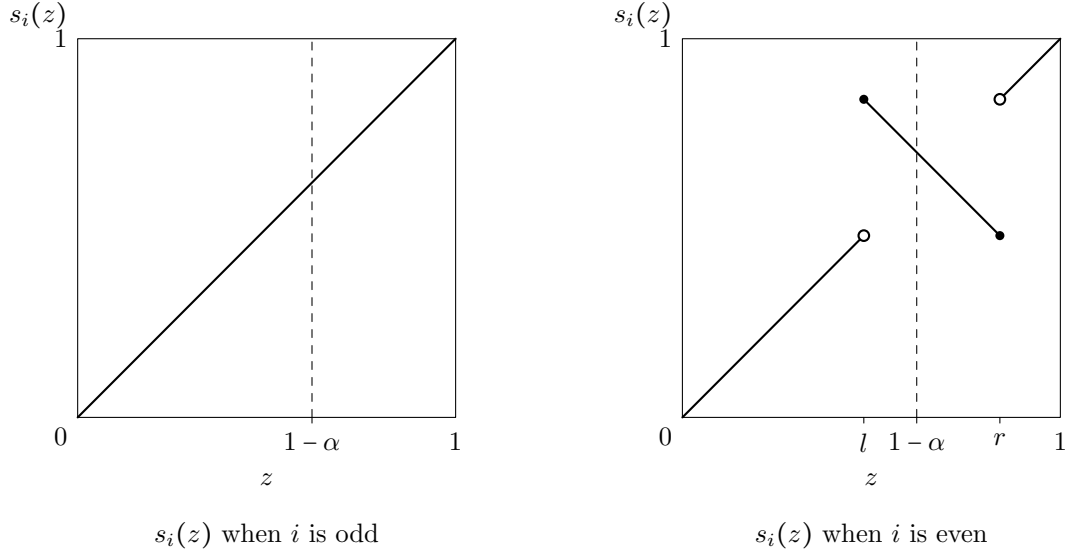
 
\centering
\vspace{-.5cm}
\include{tikz-plot-1}
\vspace{-1cm}
\caption{An illustration of the score functions $s_i(z)$ constructed in the proof of Proposition~\ref{prop:tightness}, for some choice of $l$ and $r$.}
\label{fig:tightness-scores}
\end{figure}

Without loss of generality,\footnote{
We use the fact that there is a measure-isomorphism between an atom-less Borel measurable space to the uniform Lebesgue measure space on $[0,1]$. See \cite[Thm.~17.41]{kechris2012classical}.
} we assume $P=\mathsf{Unif}(0,1)$. For $i\geq 1$, we construct the sequence of score functions $s_i(\cdot; Z_{<i}) = s_i(\cdot)$ independent of the data, as follows. Fix some $l < 1-\alpha < r$ and define
\[s_i(z) = \begin{cases} z, & \textnormal{if $i$ is odd},\\ z\cdot\ind\{z\not\in[l,r]\} + (l+r-z)\cdot\ind\{z\in[l,r]\}, &\textnormal{if $i$ is even}.\end{cases} \]
(See Figure~\ref{fig:tightness-scores} for an illustration.) Define also the average of the two cases,
\[\bar{s}(z) = \frac{s_1(z)+s_2(z)}{2} = z\cdot\ind\{z\not\in[l,r]\} + \frac{l+r}{2}\cdot\ind\{z\in[l,r]\}. \]

By definition of these score functions, along with the fact that $Z_i\sim\mathsf{Unif}(0,1)$, we have $s_i(Z_i)\sim\mathsf{Unif}(0,1)$ regardless of whether $i$ is odd or even. Therefore, we have
\[\P\prn{s_i(z)>s_i(Z_i)} = s_i(z),\]
for any fixed $z\in[0,1]$, and any $i$.
By the Dvoretzky--Kiefer--Wolfowitz inequality, since the $Z_i$'s are i.i.d., we therefore have
\[ \sup_{z\in[0,1]} \left|\frac{1}{n}\sum_{i=1}^n \ind \brc{s_i(z) > s_i(Z_i)} - \bar{s}(z)\right|\to 0,\]
almost surely. 

Next we split into cases. First fix any $\nu\in(0,\alpha)$, so that our target coverage level is $1-2\alpha+\nu\in (1-2\alpha,1-\alpha)$. Define $l = 1 - 2\alpha +\nu$ and $r = 1$.
Note that $\bar{s}(z)$ satisfies
\[\textnormal{$\bar{s}(z)\geq \frac{l+r}{2}$ if $z\geq l$, and $\bar{s}(z)\leq l$ if $z<l$},\]
where we note that $\frac{l+r}{2}>1-\alpha$ while $l < 1-\alpha$.
Consequently, 
since $z\in\mc{C}_n$ if and only if $ \sum_{i=1}^n \ind \brc{s_i(z) > s_i(Z_i)} <(1-\alpha)(n+1)$, we therefore see that
\[\mc{C}_n = [0,l)\textnormal{ for all sufficiently large $n$, almost surely}.\]
Therefore,
\[\lim_{n\to\infty}\P\prn{Z_{n+1}\in\mc{C}_n} = l = 1-2\alpha + \nu.\]

Next, suppose instead that $\nu\in(\alpha,2\alpha)$, so that our target coverage level is $1-2\alpha+\nu\in(1-\alpha,1)$.
In this case, set $l = 1-2\alpha-\epsilon$ (for some $\epsilon\in (0,1-2\alpha)$), and $r=1-2\alpha+\nu$.
In this case, $\bar{s}(z)$ satisfies
\[\textnormal{$\bar{s}(z)\leq \frac{l+r}{2}$ if $z\leq r$, and $\bar{s}(z)\geq r$ if $z>r$},\]
where we note that $\frac{l+r}{2} = 1-2\alpha + \frac{\nu-\epsilon}{2} < 1-\alpha$, while $r>1-\alpha$. Consequently, 
since $z\in\mc{C}_n$ if and only if $ \sum_{i=1}^n \ind \brc{s_i(z) > s_i(Z_i)} <(1-\alpha)(n+1)$, we therefore see that
\[\mc{C}_n = [0,r]\textnormal{ for all sufficiently large $n$, almost surely}.\]
Therefore,
\[\lim_{n\to\infty}\P\prn{Z_{n+1}\in\mc{C}_n} = r = 1-2\alpha + \nu.\]

Finally, we consider the two remaining cases. For $\nu = 2\alpha$ (i.e., coverage $1-2\alpha+\nu=1$), 
we can simply consider the example $s_i(z) \equiv 1$ for all $i\geq 1$, so that $\mc{C}_n=\mc{Z}$ almost surely, and coverage is therefore $\P\prn{Z_{n+1}\in\mc{C}_n}=1$ for all $n$. And, for $\nu=\alpha$ (i.e., coverage $1-2\alpha+\nu=1-\alpha$), we can take $s_i(z) = z$ for all $i\geq 1$, in which case the prediction sets $\mc{C}_n$ converge to $[0,1-\alpha]$, with coverage $\P\prn{Z_{n+1}\in\mc{C}_n}\to 1-\alpha$.

%% file: tikz-plot-1.tex
\begin{tikzpicture}[
    x=5cm,
    y=5cm,
    line cap=round,
    line join=round,
    every node/.style={font=\small}
]

\pgfmathsetmacro{\A}{0.62}      
\pgfmathsetmacro{\L}{0.48}     
\pgfmathsetmacro{\R}{0.84}     

\begin{scope}

    \draw (0,0) rectangle (1,1);

    \draw[thick] (0,0) -- (1,1);

    \draw[dashed] (\A,0) -- (\A,1);

    \node[above left] at (0,1) {$s_i(z)$};
    \node[below] at (0.5,-0.12) {$z$};

    \node[below left] at (0,0) {$0$};
    \node[left] at (0,1) {$1$};


    \draw (\A,0) -- ++(0,-0.015);
    \node[below] at (\A,-0.015) {$1-\alpha$};

    \node[below] at (1,-0.015) {$1$};

    \node[below] at (0.5,-0.25)
        {$s_i(z)$ when $i$ is odd};

\end{scope}

\begin{scope}[xshift=8cm]

    \draw (0,0) rectangle (1,1);

    \draw[thick]
        (0,0) -- (\L,\L);

    \draw[thick]
        (\L,\R) -- (\R,\L);

    \draw[thick]
        (\R,\R) -- (1,1);

    \draw[dashed]
        (\A,0) -- (\A,1);

    \fill (\L,\R) circle[radius=0.012];
    \fill (\R,\L) circle[radius=0.012];

    \draw[fill=white, thick]
        (\L,\L) circle[radius=0.014];

    \draw[fill=white, thick]
        (\R,\R) circle[radius=0.014];

    \node[above left] at (0,1) {$s_i(z)$};
    \node[below] at (0.5,-0.12) {$z$};

    \node[below left] at (0,0) {$0$};
    \node[left] at (0,1) {$1$};


    \draw (\L,0) -- ++(0,-0.015);
    \node[below] at (\L,-0.015) {$l$};

    \draw (\A,0) -- ++(0,-0.015);
    \node[below] at (\A,-0.015) {$1-\alpha$};

    \draw (\R,0) -- ++(0,-0.015);
    \node[below] at (\R,-0.015) {$r$};

    \node[below] at (1,-0.015) {$1$};

    \node[below] at (0.5,-0.25)
        {$s_i(z)$ when $i$ is even};

\end{scope}

\end{tikzpicture}

%% file: proof-stability.tex
\section{Proofs for Section~\ref{sec:stability}} \label{proof:stability}

\subsection{Proof of Theorem~\ref{thm:stability-score}}
Let $\gamma >0$. By Markov's inequality, it follows from Assumption~\ref{assmp:score-stability} that
\begin{align*}
     &\P \prn{\left|\sum_{i=m}^n \ind\{s_{i}(Z_{n+1}; Z_{<i}) > s_{i}(Z_i; Z_{<i})\} -  \sum_{i=m}^n  \ind\{s_\star(Z_{n+1};Z_{<m}) > s_\star(Z_i;Z_{<m})\} \right| \geq (n-m+1)\gamma} \nonumber \\
     & \leq \frac{1}{(n-m+1)\gamma } \sum_{i=m}^n \P\prn{\ind \brc{s_i(Z_{n+1}; Z_{<i}) > s_i(Z_i; Z_{<i})} \neq \ind \brc{s_\star(Z_{n+1};Z_{<m}) > s_\star(Z_i;Z_{<m})}} \\
     & =\frac{1}{(n-m+1)\gamma } \sum_{i=m}^n \P\prn{\ind \brc{s_i(Z; Z_{<i}) > s_i(Z'; Z_{<i})} \neq \ind \brc{s_\star(Z;Z_{<m}) > s_\star(Z';Z_{<m})}}
     \leq \frac{\nu}{\gamma}.
\end{align*}
Since $Z_1, \dots, Z_n, Z_{n+1}$ are i.i.d., the scores $s_\star(Z_m;Z_{<m}),\dots,s_\star(Z_{n+1};Z_{<m})$ are i.i.d.\ conditionally on $Z_{<m}$, and therefore exchangeable, we have that
\[\frac{1+\sum_{i=m}^n  \ind\{s_\star(Z_{n+1};Z_{<m}) \leq s_\star(Z_i;Z_{<m})\} }{n-m+2} \stg \mathsf{Unif}(0,1)\]
(note that the quantity on the left-hand side is simply the split conformal p-value, for the trained score function $s_\star(\cdot;Z_{<m})$, with calibration set $Z_m,\dots,Z_n$ and test point $Z_{n+1}$). Consequently,
\begin{align*}
    \P \prn{\sum_{i=m}^n  \ind\{s_\star(Z_{n+1};Z_{<m}) > s_\star(Z_i;Z_{<m})\} <  \beta(n-m+2) } \geq \beta.
\end{align*}
for any $\beta\in[0,1]$.

Combining the above calculations (with the appropriate choice of $\beta$) yields
\begin{align*}
     \P(Z_{n+1} \in \mc{C}_{n}) &= \P\prn{\sum_{i=1}^n \ind\{s_{i}(Z_{n+1}; Z_{<i}) > s_{i}(Z_i; Z_{<i})\} < (1-\alpha)(n+1)}\\
     &\geq \P\prn{\sum_{i=m}^n \ind\{s_{i}(Z_{n+1}; Z_{<i}) > s_{i}(Z_i; Z_{<i})\} < (1-\alpha)(n+1) - (m-1)}\\
     & \geq \P \prn{\sum_{i=m}^n  \ind\{s_\star(Z_{n+1}) > s_\star(Z_i)\} < (1-\alpha)(n+1) - (n-m+1)\gamma - (m-1)} - \frac{\nu}{\gamma} \nonumber \\
     & \geq \frac{(1-\alpha)(n+1) - (n-m+1)\gamma - (m-1)}{n-m+2} -  \frac{\nu}{\gamma} \geq  1-\alpha\cdot \frac{n+1}{n-m+2} -\gamma - \frac{\nu}{\gamma}.
\end{align*}
Choosing $\gamma = \sqrt{\nu}$ proves the lower bound. For the upper bound, suppose that $s_\star$ has no ties. Conditional on $Z_{<m}$, the comparison count $\sum_{i=m}^n\ind\{s_\star(Z_{n+1};Z_{<m})>s_\star(Z_i;Z_{<m})\}$ is uniform on $\{0,\dots,n-m+1\}$. Combining this fact with the bounds above yields, for any $\gamma>0$,
\begin{align*}
    \P(Z_{n+1}\in\mc{C}_n)
    &\leq \P\prn{\sum_{i=m}^n\ind\{s_i(Z_{n+1};Z_{<i})>s_i(Z_i;Z_{<i})\}<(1-\alpha)(n+1)}\\
    &\leq \P\prn{\sum_{i=m}^n\ind\{s_\star(Z_{n+1};Z_{<m})>s_\star(Z_i;Z_{<m})\}<(1-\alpha)(n+1)+(n-m+1)\gamma}+\frac{\nu}{\gamma}\\
    &\leq \frac{(1-\alpha)(n+1)+(n-m+1)\gamma+1}{n-m+2}+\frac{\nu}{\gamma}\\
    &\leq 1-\alpha+\frac{(1-\alpha)(m-1)+1}{n-m+2}+\gamma+\frac{\nu}{\gamma}.
\end{align*}
The same choice of $\gamma = \sqrt{\nu}$ proves the upper bound. This completes the proof.

\subsection{Proof of Proposition~\ref{prop:convergence-Lq-implications}}
We will prove the result for $q<\infty$; the result for $q=\infty$ can be derived similarly.

For any $i\in\{m,\dots,n\}$, conditional on $Z_{<m}$, define
\[\wt{\nu}_{q,i}(Z_{<m}) = \norm{s_i(\cdot;Z_{<i}) - s_\star(\cdot;Z_{<m})}_{L_q},\]
and
\[\nu_i(Z_{<m}) = \P\prn{\ind \brc{s_i(Z; Z_{<i}) > s_i(Z'; Z_{<i})} \neq \ind \brc{s_\star(Z; Z_{<m}) > s_\star(Z'; Z_{<m})} \mid Z_{<m}}.\]
Then, by Assumption~\ref{assmp:score-Lq-stability}, we have
\[\wt{\nu}_q^q \geq \frac{1}{n-m+1} \sum_{i=m}^n\E\brk{\big(\wt{\nu}_{q,i}(Z_{<m})\big)^q},\]
while Assumption~\ref{assmp:score-stability} holds with
\[\nu = \frac{1}{n-m+1} \sum_{i=m}^n\E\brk{\nu_i(Z_{<m})} .\]
From this point on, then, our goal will be to bound $\nu_i(Z_{<m})$ in terms of $\wt{\nu}_{q,i}(Z_{<m})$.

Fix any $i\in\{m,\dots,n\}$. For convenience we will write
\[V = s_i(Z; Z_{<i}), \ V' = s_i(Z'; Z_{<i}), \ W = s_\star(Z; Z_{<m}), \ W' = s_\star(Z'; Z_{<m}).\]
Note that $(V,W)\eqd (V',W')$, and we have
\[\wt{\nu}_{q,i}(Z_{<m}) = \norm{V-W}_{L_q} = \E\brk{|V-W|^q}^{1/q},\quad \nu_i(Z_{<m}) = \P\prn{\ind\brc{V>V'}\neq \ind\brc{W>W'}},\]
where from this point on we condition on $Z_{<m}$ implicitly. We can observe that, for any fixed $\delta>0$,
\[\textnormal{If $\ind\brc{V>V'}\neq \ind\brc{W>W'}$ then either $|W-W'| \leq \delta$ or $|V-W|\geq \delta/2$ or $|V'-W'|\geq \delta/2$.}\]
We also have
\[\P\prn{|V-W|\geq \delta/2} \leq \frac{\E\brk{|V-W|^q}}{(\delta/2)^q} \leq \left(\frac{2\wt{\nu}_{q,i}(Z_{<m})}{\delta}\right)^q,\]
and the same bound holds for $\P\prn{|V'-W'|\geq \delta/2}$.
Moreover, since $W,W'$ are i.i.d.\ and have density bounded by $B$ (conditional on $Z_{<m}$), 
\[\P\prn{|W-W'|\leq \delta} = \E\brk{\P\prn{W \in [W'-\delta,W'+\delta]\mid W'}} \leq \E\brk{B\cdot 2\delta} = 2\delta B.  \]
Combining everything, then,
\[\nu_i(Z_{<m}) =  \P\prn{\ind\brc{V>V'}\neq \ind\brc{W>W'}} \leq 2\left(\frac{2\wt{\nu}_{q,i}(Z_{<m})}{\delta}\right)^q +  2\delta B.\]

Returning to our work above, we then see that Assumption~\ref{assmp:score-stability} holds with
\begin{multline*}\nu =\frac{1}{n-m+1} \sum_{i=m}^n\E\brk{\nu_i(Z_{<m})} 
\leq \frac{1}{n-m+1} \sum_{i=m}^n\E\brk{ 2\left(\frac{2\wt{\nu}_{q,i}(Z_{<m})}{\delta}\right)^q +  2\delta B}\\
= 2\delta B + \frac{2^{q+1}}{\delta^q} \cdot\frac{1}{n-m+1}\sum_{i=m}^n \E\brk{\big(\wt{\nu}_{q,i}(Z_{<m})\big)^q} \leq 2\delta B + \frac{2^{q+1}}{\delta^q} \cdot \wt{\nu}_q^q.\end{multline*}
Choosing $\delta = \frac{(2\wt{\nu}_q)^{q/(q+1)}}{B^{1/(q+1)}}$ and noticing that $\nu \le 1$ always hold completes the proof.

%% file: appendix_additional_results.tex
\section{Additional results and extensions}\label{appendix_additional_results}

\subsection{Coverage under distribution drift}
Throughout this paper, we have assumed that the data points $Z_1,\dots,Z_{n+1}$ are i.i.d.\ from a common distribution $P$ (or, more generally, exchangeable). In this section, we now turn to a relaxation of this condition. Suppose that the data points are independent, with $Z_i \sim P_i$ for each $i$. At a high level, our goal is to show that rolling-CP is robust to mild distribution drift---if most distributions $P_i$ are similar to the distribution $P_{n+1}$ of the test point, then rolling-CP should not lose much coverage. (Analogous results for full conformal prediction can be found in the work of \citep{barber2023conformal} on ``nonexchangeable conformal prediction''.)

\begin{theorem}\label{thm:distribution-drift}
    Let $Z_1,\dots,Z_{n+1}$ be independent, with $Z_i\sim P_i$ for each $i\in[n+1]$. Let $\alpha\in(0,1)$. For any sequence of score functions $s_i:\mc{Z}\times\mc{Z}^{i-1}\to\R$,
the rolling-CP prediction set $\mc{C}_n$ defined in~\eqref{eqn:iterative-conformal-set} satisfies
\begin{align}
\P\prn{Z_{n+1}\in\mc{C}_n}\geq 1-2\alpha - \frac{2}{n}\sum_{i=1}^n \mathsf{d}_{\mathsf{TV}}(P_i,P_{n+1}),
\end{align}
where $\mathsf{d}_{\mathsf{TV}}$ denotes the total variation distance between distributions.
\end{theorem}
\noindent In other words, if most distributions $P_i$ are close to the test point distribution $P_{n+1}$ (with respect to total variation distance), then the coverage result is nearly the same as in the exchangeable or i.i.d.\ setting.
\paragraph{Proof of Theorem~\ref{thm:distribution-drift}.}
    The proof follows similar strategies as the proof of Theorem~\ref{thm:exchangeable-marginal} for the i.i.d.\ case (presented in Section~\ref{sec:marginal-coverage}). 

    Let $\mc{A}_i = \brc{s_i(Z_{n+1};Z_{<i})>s_i(Z_i;Z_{<i})}$, and let $p^*_i = \P\prn{\mc{A}_i\mid Z_{<i}} = \P\prn{\mc{A}_i\mid Z_{<n+1}}$. Define also
    \[p^{**}_i = \P\prn{s_i(Z_i';Z_{<i})>s_i(Z_i;Z_{<i}) \mid Z_{<i+1}},\]
    where $Z'_i\sim P_i$ is independent of the data. 
    Then
    \begin{multline*}|p^*_i - p^{**}_i| = \left|\P\prn{s_i(Z_{n+1};Z_{<i})>s_i(Z_i;Z_{<i}) \mid Z_{<n+1}} - \P\prn{s_i(Z'_i;Z_{<i})>s_i(Z_i;Z_{<i}) \mid Z_{<n+1}}\right|\\\leq \mathsf{d}_{\mathsf{TV}}(Z'_i,Z_{n+1}) = \mathsf{d}_{\mathsf{TV}}(P_i,P_{n+1}),\end{multline*}
    almost surely. To account for the possibility of ties, next define
    \[p_i = U_i \cdot \P\prn{s_i(Z'_i;Z_{<i})=s_i(Z_i;Z_{<i}) \mid Z_{<i+1}} + \P\prn{s_i(Z'_i;Z_{<i})>s_i(Z_i;Z_{<i}) \mid Z_{<i+1}} \]
    for $U_1,\dots,U_n\simiid\mathsf{Unif}(0,1)$ (drawn independently of the data).  Then
    we see that $p_i \mid Z_{<i}\sim \mathsf{Unif}(0,1)$ for each $i$, and consequently $p_1,\dots,p_n \simiid\mathsf{Unif}(0,1)$ (by the same argument as in the proof of Theorem~\ref{thm:exchangeable-marginal} for the i.i.d.\ case given in Section~\ref{sec:marginal-coverage}, and as in the proof of Theorem~\ref{thm:training-conditional}). Moreover, by construction, $p_i \geq p^{**}_i$ almost surely, and so $p^*_i\leq p_i + \mathsf{d}_{\mathsf{TV}}(P_i,P_{n+1})$, almost surely.

        Our next step is a lemma, similar to the results of Lemmas~\ref{lem:pigeonhole} and~\ref{lem:pigeonhole_conditional}:
    \begin{lemma}\label{lem:pigeonhole_drift}
    Let $\mc{A}_1,\dots,\mc{A}_n$ be any events. Suppose that there exists some random variable $W$, constants $c_1,\dots,c_n\geq 0$, and random variables $p_1,\dots,p_n$, such that 
    \[\P\prn{\mc{A}_i \mid W}\leq p_i + c_i\textnormal{ almost surely for all $i\in[n]$, and }p_1,\dots,p_n\simiid \mathsf{Unif}(0,1).\]
    Then for any $\Delta \geq 0$,
    \[\P\prn{\sum_{i=1}^n \ind_{\mc{A}_i} \geq  n-\Delta } \leq \frac{2\Delta+1}{n+1} + \frac{2}{n}\sum_{i=1}^n c_i.\]
\end{lemma}
We the apply this lemma with $\Delta = \alpha(n+1) - 1$, $W=Z_{<n+1}$, and $c_i = \mathsf{d}_{\mathsf{TV}}(P_i,P_{n+1})$. This completes the proof, since as before, $Z_{n+1}\in\mc{C}_n$ if and only if $\sum_{i=1}^n\ind_{\mc{A}_i} < (1-\alpha)(n+1)$.
\endproof

\subsection{Sequential inference via merging: a comparison}
As mentioned in Section~\ref{sec:related-work}, we can think of rolling-CP as an approach that merges information from the $n$ trained score functions, $s_i(\cdot;Z_{<i})$ for $i\in[n]$, to produce a single prediction set. We may therefore ask how this compares to other methods for merging outputs from multiple predictive inference procedures.

In particular, for each value $i$, we could consider the split-CP prediction set obtained with the $i$-th trained score function: recalling~\eqref{eqn:splitCP} (and updating notation for data indexing as needed), this prediction set is given by
\[\mc{C}_{\textnormal{split-CP},i} = \left\{z : \sum_{j=i}^n \ind\brc{s_i(z;Z_{<i})>s_i(Z_j;Z_{<i})} < (1-\alpha)(n-i+2)\right\}.\]
This can equivalently be expressed as
\[\mc{C}_{\textnormal{split-CP},i} = \left\{z : p_{\textnormal{split-CP},i}(z) > \alpha\right\} \textnormal{ where } p_{\textnormal{split-CP},i}(z) = \frac{1 + \sum_{j=i}^n \ind\brc{s_i(z;Z_{<i})\leq s_i(Z_j;Z_{<i})}}{n-i+2} .\]
This quantity $p_{\textnormal{split-CP},i}(z)$ is the split conformal p-value \citep{vovk2005algorithmic}, and is a valid p-value at the test point $Z_{n+1}$---that is, for exchangeable data, $p_{\textnormal{split-CP},i}(Z_{n+1})\stg \mathsf{Unif}(0,1)$.

Next, fix any weights $w_1,\dots,w_n\geq 0$ with $\sum_{i=1}^n w_i=1$. Let
\[\bar{p}_{\textnormal{split-CP}} = \sum_{i=1}^n w_i \cdot p_{\textnormal{split-CP},i}(z)\]
be the (weighted) average of the split-CP p-values. Defining
$\bar{\mc{C}}_{\textnormal{split-CP}} = \{z : \bar{p}_{\textnormal{split-CP}}>\alpha\}$,
we then have the guarantee
\[\P\prn{Z_{n+1}\in \bar{\mc{C}}_{\textnormal{split-CP}} }\geq 1-2\alpha,\]
because any average of p-values is a p*-value and consequently offers validity up to a factor-of-two  \citep{wang2024testing}.

This appears to be a very similar construction, and the same guarantee, as rolling-CP. In particular, each method offers the ability to train the model on all data points arriving sequentially in the data stream. However, there is a key difference: while rolling-CP places equal weight on the comparison between the test point and each $Z_i$, the merged split-CP prediction set has a lower effective sample size because it places unequal weights on different training points: we can calculate
\[\bar{p}_{\textnormal{split-CP}} = \sum_{i=1}^n w_i \cdot p_{\textnormal{split-CP},i}(z) = \sum_{i=1}^n w_i \cdot \frac{1 + \sum_{j=i}^n \ind\brc{s_i(z;Z_{<i})\leq s_i(Z_j;Z_{<i})}}{n-i+2}, \]
so the total weight placed on the comparison of $z$ with $Z_j$ is given by
\[\sum_{i=1}^j \frac{w_i}{n-i+2},\]
which is increasing with $j$. In other words, the merged split-CP method places substantially more weight on later data points, and may then lose accuracy due to lower effective sample size. (For rolling-CP, which places equal weight on each data point $Z_1,\dots,Z_n$, there is no immediately apparent way to prove its coverage result using an analogous average-of-p-values style argument.)

\subsection{Training-conditional coverage under stability}
Finally, we present an additional result in the setting of stability, extending the findings of Section~\ref{sec:stability}: 
we show that training-conditional coverage is $\gtrapprox 1-\alpha$, under the stability condition of Assumption~\ref{assmp:score-stability}. This can be viewed as a stronger form of Theorem~\ref{thm:training-conditional} (which establishes training-conditional coverage $\gtrapprox 1-2\alpha$, without a stability condition).

\begin{theorem}[Training-conditional coverage under score comparison stability] \label{thm:stability-score-training-conditional} 
Let $Z_1, \cdots, Z_n, Z_{n+1} \simiid P$ and let $\alpha \in (0,1)$. Define the rolling-CP prediction set $\mc{C}_n$ as in~\eqref{eqn:iterative-conformal-set}, and suppose Assumption~\ref{assmp:score-stability} holds for some $m\in[n]$ and some $\nu\in[0,1]$. Let
    \[\alpha_P(Z_{< n+1}) = \P(Z_{n+1} \notin \mc{C}_{n} \mid Z_1, \dots, Z_n)\]
    as the training-conditional miscoverage rate of $\mc{C}_n$. Then, for any $\delta \in (0,1)$,
    \begin{equation*}
        \P\prn{\alpha_P(Z_{<n+1}) \leq \alpha \cdot \frac{n+1}{n-m+2} + \sqrt{\frac{\log(1/\delta)}{2(n-m+2)}}+2 \sqrt[3]{\nu}} \geq  1 - \delta - \sqrt[3]{\nu}.
    \end{equation*}
\end{theorem}

\paragraph{Proof of Theorem~\ref{thm:stability-score-training-conditional}.}
The proof combines ideas from the arguments behind Theorem~\ref{thm:training-conditional} (training-conditional coverage without stability), and Theorem~\ref{thm:stability-score} (marginal coverage, with stability).

Fix any $\gamma>0$. As in the proof of Theorem~\ref{thm:stability-score}, it holds that
\[\P \prn{\left|\sum_{i=m}^n \ind\{s_{i}(Z_{n+1}; Z_{<i}) > s_{i}(Z_i; Z_{<i})\} -  \sum_{i=m}^n  \ind\{s_\star(Z_{n+1};Z_{<m}) > s_\star(Z_i;Z_{<m})\} \right| \geq (n-m+1)\gamma }\leq \frac{\nu}{\gamma}.\]
Define a random variable
\[D = \sum_{i=m}^n \ind\{s_{i}(Z_{n+1}; Z_{<i}) > s_{i}(Z_i; Z_{<i})\} - \sum_{i=m}^n  \ind\{s_\star(Z_{n+1};Z_{<m}) > s_\star(Z_i;Z_{<m})\},\]
so that the bound above implies
\[\P\prn{D \geq (n-m+1)\gamma}\leq \frac{\nu}{\gamma}.\]
Since $Z_{n+1}\in\mc{C}_n$ if and only if $\sum_{i=1}^n \ind\{s_{i}(Z_{n+1}; Z_{<i}) > s_{i}(Z_i; Z_{<i})\} <(1-\alpha)(n+1)$, it holds that
\[\textnormal{If $Z_{n+1}\not\in\mc{C}_n$ then $\sum_{i=m}^n  \ind\{s_\star(Z_{n+1};Z_{<m}) > s_\star(Z_i;Z_{<m})\} + D \geq (1-\alpha)(n+1)-(m-1)$.}\]
Therefore,
\begin{multline*}\alpha_P(Z_{<n+1}) \leq \P\prn{\sum_{i=m}^n  \ind\{s_\star(Z_{n+1};Z_{<m}) > s_\star(Z_i;Z_{<m})\} + D \geq (1-\alpha)(n+1)-(m-1) \,\Big|\, Z_{<n+1}}\\
\leq \P\prn{\sum_{i=m}^n  \ind\{s_\star(Z_{n+1};Z_{<m}) > s_\star(Z_i;Z_{<m})\}\geq (1-\alpha)(n+1) -(m-1)- (n-m+1)\gamma  \,\Big|\, Z_{<n+1}}\\
+ \P\prn{D \geq  (n-m+1)\gamma  \,\Big|\, Z_{<n+1}}.
\end{multline*}
By Markov's inequality, for any $c> 0$,
\[\P\prn{ \P\prn{D \geq  (n-m+1)\gamma  \,\Big|\, Z_{<n+1}} \geq c} \leq \frac{\E\brk{\P\prn{D \geq (n-m+1)\gamma \,\Big|\, Z_{<n+1}}}}{c} = \frac{\P\prn{D \geq  (n-m+1)\gamma}}{c} \leq \frac{\nu}{c\gamma}.\]
Next note that $s_\star(Z_m;Z_{<m}),\dots,s_\star(Z_{n+1};Z_{<m})$ are i.i.d.\ conditional on $Z_{<m}$. Denote $S_i = s_\star(Z_i;Z_{<m})$ for each $i = m,\cdots, n+1$. Let $S_{(1)}\leq \dots \leq S_{(n-m+1)}$ be the order statistics of $S_m, \cdots, S_n$. 
Then
\[\sum_{i=m}^n  \ind\{s_\star(Z_{n+1};Z_{<m}) > s_\star(Z_i;Z_{<m})\}\geq (1-\alpha)(n+1) -(m-1)- (n-m+1)\gamma \ \Longleftrightarrow \ S_{n+1} > S_{(k)}\] where $k = \lceil (1-\alpha)(n+1) -(m-1)- (n-m+1)\gamma\rceil$. Writing $F$ as the CDF of $s_\star(Z;Z_{<m})$ (for $Z\sim P$, conditional on $Z_{<m}$), we then have
\[\P\prn{S_{n+1} > S_{(k)}\mid Z_{<n+1}} = 1 - F(S_{(k)}).\]
Since $S_{(k)}$ is an order statistic of $S_m,\dots,S_n$, which are i.i.d.\ draws from the CDF $F$ (after conditioning on $Z_{<m}$), 
\[F(S_{(k)})\stg U_{(k)}, \]
where $U_{(k)}$ is the $k$th order statistic of $n-m+1$ i.i.d.\ $\mathsf{Unif}(0,1)$ random variables, which satisfies $U_{(k)} \sim\mathsf{Beta}(k,n-m+2-k)$. Therefore, for $c'>0$,
\[\P\prn{\P\prn{S_{n+1} > S_{(k)}\mid Z_{<n+1}}  \geq \frac{n-m+2-k}{n-m+2} + c'}
\leq \P\prn{ U_{(k)} \leq \frac{k}{n-m+2} - c'} \leq e^{-2(n-m+2)c'{}^2},\]
by subgaussianity of the Beta distribution \citep{elder2016bayesian}. 

Combining everything,
\[\P\prn{\alpha_P(Z_{<n+1}) \leq \frac{n-m+2-k}{n-m+2} + c' + c} \geq  1 - \frac{\nu}{c\gamma} - e^{-2(n-m+2)c'{}^2}.\]
Plugging in our choice of $k$, and simplifying,
\[\P\prn{\alpha_P(Z_{<n+1}) \leq \frac{\alpha(n+1) +(n-m+1)\gamma}{n-m+2} + c' + c} \geq  1 - \frac{\nu}{c\gamma} - e^{-2(n-m+2) c'{}^2}.\]
Choosing $\gamma = c = \sqrt[3]{\nu}$ and $c' = \sqrt{\frac{\log(1/\delta)}{2(n-m+2)}}$, 
this simplifies to
\[\P\prn{\alpha_P(Z_{<n+1}) \leq \alpha \cdot \frac{n+1}{n-m+2} + \sqrt{\frac{\log(1/\delta)}{2(n-m+2)}} + 2\sqrt[3]{\nu}} \geq  1 -\delta- \sqrt[3]{\nu} ,\]
which completes the proof.
\endproof

%% file: appendix_lemmas.tex
\section{Proofs of lemmas}\label{appendix_lemmas}

\subsection{Proof of Lemma~\ref{lem:pigeonhole_exch}}
First for intuition, we recall the simpler form of the proof of Lemma~\ref{lem:pigeonhole}. In that result, we had $p_i = \P\prn{\mc{A}_i\mid W}$ for a single shared random variable $W$, and we were then able to apply the pigeonhole principle to reduce to a set $p_{(1)},\dots,p_{(k)}$, the $k$ smallest values among $p_1,\dots,p_n$. In this new setting, however, we cannot follow the same argument. For instance, in order to determine whether $p_i$ belongs to this set (i.e., whether $p_i$ is one of the $k$ smallest values), we must observe \emph{all} the random variables $p_1,\dots,p_n$ for comparison. In the i.i.d.\ setting, since $p_i = \P\prn{\mc{A}_i\mid Z_{<n+1}}$ while the random variables $p_1,\dots,p_n$ are $Z_{<n+1}$-measurable, we are allowed to observe $p_1,\dots,p_n$---that is, it still holds that $p_i$ is the conditional probability of the event $\mc{A}_i$, even if we observe $p_1,\dots,p_n$. In contrast, in the exchangeable case, we have $p_i = \P\prn{\mc{A}_i\mid \mc{F}_i}$, but $p_{i+1},\dots,p_n$ are \emph{not} $\mc{F}_i$-measurable; after observing the random variables $p_1,\dots,p_n$, it is no longer the case that $p_i$ is the conditional probability of $\mc{A}_i$.

Consequently, we will apply the pigeonhole principle in a more subtle way. We will define random variables $B_1,\dots,B_n\in\{0,1\}$, adapted to the filtration, so that $B_i = 1$ indicates that $p_i$ is selected for the pigeonhole argument. Fixing an integer $k\in\{0,\dots,n\}$, define
\[B_1 = \ind\left\{p_1 \leq \frac{k}{n}, \ 0 < k\right\},\]
and then inductively for each $i=2,\dots,n$,
\[B_i = \ind\left\{p_i\leq \frac{k-\sum_{j=1}^{i-1}B_j}{n-i+1}, \ \sum_{j=1}^{i-1}B_j <k\right\}.\]
Note that $p_j\in\mc{F}_j\subseteq\mc{F}_i$ for all $j\leq i$, by construction, and so $B_i \in\mc{F}_i$. 

Moreover, by construction, we must have $\sum_{i=1}^n B_i = k$, almost surely.
Therefore, by the pigeonhole principle,
\[\textnormal{If $\sum_{i=1}^n \ind_{\mc{A}_i} \geq n-\Delta$ then $\sum_{i\in[n], B_i=1}\ind_{\mc{A}_i} \geq k-\Delta$.}\]
Then,
\[\P\prn{\sum_{i=1}^n \ind_{\mc{A}_i}\geq n-\Delta}\leq \P\prn{\sum_{i\in[n], B_i=1}\ind_{\mc{A}_i} \geq k-\Delta} \leq \frac{\E\brk{\sum_{i\in[n], B_i=1}\ind_{\mc{A}_i}}}{k-\lfloor\Delta\rfloor} = \frac{\sum_{i=1}^n \E\brk{\ind_{\mc{A}_i} B_i}}{k-\lfloor\Delta\rfloor},\]
as long as $k\geq \lfloor\Delta\rfloor$,
where as before we apply Markov's inequality, together with the fact that the sum $\sum_{i\in[n], B_i=1}\ind_{\mc{A}_i}$ is integer-valued (so that we can replace $\Delta$ with $\lfloor \Delta\rfloor$ in the denominator). We can also calculate
\[\E\brk{\ind_{\mc{A}_i} B_i} = \E\brk{\P\prn{\mc{A}_i\mid\mc{F}_i}\cdot B_i} = \E\brk{p_i B_i},\]
since $B_i\in\mc{F}_i$ while $p_i = \P\prn{\mc{A}_i\mid\mc{F}_i}$. Finally, it holds that $\sum_{i=1}^n \E\brk{p_i B_i} = \frac{k(k+1)}{2(n+1)}$, by the lemma below; choosing $k=2\lfloor \Delta \rfloor +1$ completes the proof.

\begin{lemma} \label{lem:V-n-k-solution}
For each $n\geq 1$ and $k\in\{0,\dots,n\}$, define
random variables $p_1(n),\dots,p_n(n)$ and $B_1(n,k),\dots,B_n(n,k)$ as
\[p_i(n)\sim\mathsf{Unif}\left(\left\{0,\frac{1}{n-i+1},\dots,\frac{n-i}{n-i+1},1\right\}\right),\textnormal{ independently for $i=1,\dots,n$,}\]
and inductively,
\[B_1(n,k) = \ind\left\{p_1(n) \leq \frac{k}{n}, \ 0 < k\right\}, \qquad B_i(n,k) = \ind\left\{p_i(n)\leq \frac{k-\sum_{j=1}^{i-1}B_j(n,k)}{n-i+1}, \ \sum_{j=1}^{i-1}B_j(n,k)<k\right\}.\]
    Let $V(n,k) = \sum_{i=1}^n \E\brk{B_i(n,k)p_i(n)}$. Then we have the recursive relation
   \begin{align*}
        V(n,k) &  =  \frac{k(k+1)}{2n(n+1)} + \frac{k+1}{n+1} V(n-1, k-1) + \frac{n-k}{n+1} V(n-1,k),
    \end{align*}
    for $1 \leq k \leq n-1$, and the boundary conditions $V(n,0)=0, V(n, n)=n/2$. The preceding recursive formula admits the unique solution:
    \begin{align*}
        V(n,k) = \frac{k(k+1)}{2(n+1)} \qquad \textrm{for all} \quad 0 \leq k \leq n.
    \end{align*}
\end{lemma}
\paragraph{Proof of Lemma~\ref{lem:V-n-k-solution}.}
    We first notice that $p_{i-1}(n-1)$ has exactly the same distribution as $p_i(n)$. Moreover, it is easy to see conditional on $B_1(n,k)=1$, the joint distributions below are identical:
    \begin{align*}
        &(B_2(n,k), p_2(n)), \dots, (B_n(n,k), p_n(n)) \nonumber \\
        & \qquad \stackrel{d}{=} (B_1(n-1,k-1), p_1(n-1)), \dots, (B_{n-1}(n-1,k-1), p_{n-1}(n-1)).
    \end{align*}
    While on $B_1(n,k)=0$, we have similarly
    \begin{align*}
        &(B_2(n,k), p_2(n)), \dots, (B_n(n,k), p_n(n)) \nonumber \\
        & \qquad \stackrel{d}{=} (B_1(n-1,k), p_1(n-1)), \dots, (B_{n-1}(n-1,k), p_{n-1}(n-1)).
    \end{align*}
    This implies
    \begin{align*}
        V(n,k) & = \P(B_1(n,k)=1) \prn{\E \brk{p_1(n) \mid B_1(n,k)=1} +  \E \brk{\sum_{i=2}^n B_i(n,k)p_i(n) \mid B_1(n,k)=1}} \nonumber \\
        & \qquad + \P(B_1(n,k)=0) \E \brk{\sum_{i=2}^n B_i(n,k)p_i(n) \mid B_1(n,k)=0} \nonumber \\
        & \stackrel{\mathrm{(i)}}{=} \P(B_1(n,k)=1) \prn{\E \brk{p_1(n) \mid B_1(n,k)=1} +  \sum_{i=1}^{n-1}\E \brk{ B_i(n-1,k-1)p_i(n-1)}} \nonumber \\
        & \qquad + \P(B_1(n,k)=0) \sum_{i=1}^{n-1} \E \brk{ B_i(n-1,k)p_i(n-1)} \nonumber \\
        & = \P(B_1(n,k)=1) \prn{\E \brk{p_1(n) \mid B_1(n,k)=1} + V(n-1,k-1)} + \P(B_1(n,k)=0) V(n-1,k),
    \end{align*}
 in (i) we use the preceding two equalities in distribution. Substituting the explicit decision rule $B_1(n,k)$ and the distribution of $p_1(n)$, it is clear that
    \begin{align*}
        V(n,k) &= \frac{k+1}{n+1} \cdot \prn{\frac{k}{2n} + V(n-1,k-1)} + \frac{n-k}{n+1} V(n-1,k) \nonumber \\
        & = \frac{k(k+1)}{2n(n+1)} + \frac{k+1}{n+1} V(n-1, k-1) + \frac{n-k}{n+1} V(n-1,k). 
    \end{align*}
    The boundary conditions are given by $V(n,0)=0$ and $V(n,n)=n/2$. These conditions are sufficient to uniquely determine the doubly indexed recursive sequence $\{V(n,k)\}$. Indeed, one can verify that the ansatz
\begin{align*}
    V(n,k) = \frac{k(k+1)}{2(n+1)}
\end{align*}
satisfies the recursive relation exactly and thus is its unique solution.
\endproof

\subsection{Proof of Lemma~\ref{lem:dcx_comonotone}}
    First by~\eqref{eqn:dcx_equiv} we can construct a conditional distribution $P_{X\mid A}$ such that if we sample $X'\mid A \sim P_{X\mid A}$ then $\mathbb{E}[X'\mid A] \leq A$ almost surely, and marginally (i.e., marginalizing over the distribution of $A$), we have $X'\eqd X$. Let $P_{Y\mid B}$ be defined analogously for $Y$ and $B$. Define random variables $X',Y'$ by sampling
    \[(X',Y')\mid (A,B)\sim P_{X\mid A}\times P_{Y\mid B},\]
    so that we therefore have
    \[\mathbb{E}[X'\mid A,B] \leq A, \quad \mathbb{E}[Y'\mid A,B]\leq B,\]
    almost surely, and marginally $X'\eqd X$, $Y'\eqd Y$. Therefore,
    \[\mathbb{E}[X'+Y'\mid A+B]  = \mathbb{E}[\mathbb{E}[X'+Y'\mid A,B]\mid A+B] \leq \mathbb{E}[A + B\mid A+B] = A+B,\]
    which verifies that $A+B\dcx X'+Y'$ by~\eqref{eqn:dcx_equiv}. Finally, let $(X^*,Y^*)$ be a comonotone coupling for $X',Y'$ (which, equivalently, is a comonotone coupling for $X,Y$, since $X'\eqd X$ and $Y'\eqd Y$). By \citet[Theorem 5]{cote2025convex} we have $X'+Y'\cvx X^*+Y^*$ which implies $X'+Y'\dcx X^*+Y^*$. Combining these calculations completes the proof.

\subsection{Proof of Lemma~\ref{lem:pigeonhole_conditional}}
We define the same notation as in the proof of Lemma~\ref{lem:pigeonhole}: 
write $p_{(1)}\leq \dots \leq p_{(n)}$ for the order statistics of $p_1,\dots,p_n$, and let $\mc{A}_{(j)}$ denote the event associated with $p_{(j)}$. As in the proof of Lemma~\ref{lem:pigeonhole}, by the pigeonhole principle,
\[\textnormal{If $\sum_{i=1}^n \ind_{\mc{A}_i} \geq n- \Delta$ then $\sum_{j=1}^k \ind_{\mc{A}_{(j)}}\geq k-\Delta$ and therefore $\sum_{j=1}^k \ind_{\mc{A}_{(j)}}\geq k-\lfloor\Delta\rfloor$}.\]
Then applying Markov's inequality, for any positive integer $k\geq \lfloor \Delta\rfloor$ it holds that
\[\P\prn{\sum_{i=1}^n \ind_{\mc{A}_i} \geq n- \Delta\,\bigg|\, W} \leq \frac{\E\brk{\sum_{j=1}^k \ind_{\mc{A}_{(j)}}  \,\Big|\, W}}{k-\lfloor \Delta\rfloor} \leq \frac{\sum_{j=1}^k p_{(j)}}{k-\lfloor\Delta\rfloor}.\]
Next define a function $f(p_1,\dots,p_n) = \sum_{j=1}^k p_{(j)}$. By construction, for any $p_1,\dots,p_n\in[0,1]$ and $p'_i\in[0,1]$,
\[\left|f(p_1,\dots,p_n) - f(p_1,\dots,p_{i-1},p'_i,p_{i+1},\dots,p_n)\leq 1\right|,\]
and consequently by McDiarmid's inequality,
\[\P\prn{\sum_{j=1}^k p_{(j)} \geq \E\brk{\sum_{j=1}^k p_{(j)}} + \epsilon} \leq e^{-2\epsilon^2/n}.\]
By properties of the uniform distribution, $\E\brk{p_{(j)}} = \frac{j}{n+1}$, and so choosing $\epsilon = \sqrt{\frac{n\log(1/\delta)}{2}}$, we obtain
\[\P\prn{\sum_{j=1}^k p_{(j)} \geq \frac{k(k+1)}{2(n+1)} + \sqrt{\frac{n\log(1/\delta)}{2}}} \leq \delta.\]
Combining everything, we have shown that with probability $\geq 1-\delta$, it holds that 
\[\P\prn{\sum_{i=1}^n \ind_{\mc{A}_i} \geq n- \Delta\,\bigg|\, W} \leq \frac{1}{k-\lfloor\Delta\rfloor}\left(\frac{k(k+1)}{2(n+1)} +  \sqrt{\frac{n\log(1/\delta)}{2}}\right).\]
Choosing $k = 2\lfloor\Delta\rfloor +1$ completes the proof.

\subsection{Proof of Lemma~\ref{lem:pigeonhole_drift}}
We define the same notation as in the proof of Lemma~\ref{lem:pigeonhole}: 
write $p_{(1)}\leq \dots \leq p_{(n)}$ for the order statistics of $p_1,\dots,p_n$. Let $\mc{A}_{(j)}$ denote the event associated with $p_{(j)}$, and let $c_{(j)}$ be the constant associated with $p_{(j)}$ (i.e., we emphasize that $c_{(j)}$ is not an order statistic of $c_1,\dots,c_n$). 

As in the proof of Lemma~\ref{lem:pigeonhole}, by the pigeonhole principle,
\[\textnormal{If $\sum_{i=1}^n \ind_{\mc{A}_i} \geq n- \Delta$ then $\sum_{j=1}^k \ind_{\mc{A}_{(j)}}\geq k-\Delta$ and therefore $\sum_{j=1}^k \ind_{\mc{A}_{(j)}}\geq k-\lfloor\Delta\rfloor$}.\]
Then applying Markov's inequality, for any positive integer $k\geq \lfloor \Delta\rfloor$ it holds that
\[\P\prn{\sum_{i=1}^n \ind_{\mc{A}_i} \geq n- \Delta\,\bigg|\, W} \leq \frac{\E\brk{\sum_{j=1}^k \ind_{\mc{A}_{(j)}}  \,\Big|\, W}}{k-\lfloor \Delta\rfloor} \leq \frac{\sum_{j=1}^k (p_{(j)} + c_{(j)})}{k-\lfloor\Delta\rfloor}.\]
Therefore,
\[\P\prn{\sum_{i=1}^n \ind_{\mc{A}_i} \geq n- \Delta} \leq \E\brk{\frac{\sum_{j=1}^k (p_{(j)} + c_{(j)})}{k-\lfloor\Delta\rfloor}} = \frac{k(k+1)}{2(n+1)(k-\lfloor\Delta\rfloor)} + \frac{\sum_{j=1}^k \E\brk{c_{(j)}}}{k-\lfloor\Delta\rfloor},\]
since $\E\brk{p_{(j)}} = \frac{j}{n+1}$, as before. 

Next, let $\mc{E}_{i,j}$ be the event that the value $p_i$ is ranked in position $j$, i.e., $\mc{E}_{i,j} = \ind\brc{p_i = p_{(j)}}$. Recall that the $p_i$'s are continuous and so there are no ties, almost surely; therefore, the correspondence between values $p_i$ and sorted values $p_{(j)}$ is one-to-one, almost surely. Then for each $j\in[n]$, it holds that $c_{(j)} = \sum_{i=1}^n \ind\brc{\mc{E}_{i,j}}c_j$, almost surely, and so
\[\E\brk{c_{(j)}} = \sum_{i=1}^n c_i \P\prn{\mc{E}_{i,j}}. \]
Since $p_1,\dots,p_n\simiid\mathsf{Unif}(0,1)$, the ranking of this list is uniformly random, and therefore $\P\prn{\mc{E}_{i,j}}=\frac{1}{n}$ for all $i,j\in[n]$. Thus
\[\E\brk{c_{(j)}} = \frac{1}{n}\sum_{i=1}^n c_i.\]
Combining everything, we have now shown that
\[\P\prn{\sum_{i=1}^n \ind_{\mc{A}_i} \geq n- \Delta} \leq  \frac{k(k+1)}{2(n+1)(k-\lfloor\Delta\rfloor)} + \frac{\sum_{j=1}^k \left(\frac{1}{n}\sum_{i=1}^n c_i\right)}{k-\lfloor\Delta\rfloor} = \frac{k(k+1)}{2(n+1)(k-\lfloor\Delta\rfloor)} + \frac{\frac{k}{n}\sum_{i=1}^n c_i}{k-\lfloor\Delta\rfloor}.\]
Choosing $k = 2\lfloor \Delta\rfloor +1$ as before, we have completed the proof.